\documentclass [11pt]{amsart}
\usepackage {amsmath, amssymb, a4wide, epsfig, enumerate, psfrag}
\usepackage[latin1]{inputenc}
\date{\today}

 \author{Romain Dujardin}
 \title[Measurable view on parameter space]{Cubic polynomials: 
a measurable view on parameter space}


\address{Universit{\'e} Paris 7 et Institut de Math{\'e}matiques de Jussieu
         {\'E}quipe G{\'e}om{\'e}trie et Dynamique,
         Case 7012, 2  place Jussieu,
         75251 Paris Cedex 05
         France}
\email{dujardin@math.jussieu.fr}

\newcommand{\cc}{\mathbb{C}}

\newcommand{\dd}{\mathbb{D}}

\newcommand{\nn}{\mathbb{N}}
\newcommand{\e}{\varepsilon}
\newcommand{\cv}{\rightarrow}

\newcommand{\fr}{\partial}
\newcommand{\om}{\Omega}
\newcommand{\set}[1]{\left\{#1\right\}}
\newcommand{\norm}[1]{\left\Vert#1\right\Vert}
\newcommand{\abs}[1]{\left\vert#1\right\vert}
\newcommand{\cd}{\cc^2}
\newcommand{\pd}{{\mathbb{P}^2}}
\newcommand{\pu}{{\mathbb{P}^1}}
\newcommand{\rest}[1]{ \arrowvert_{#1}}

\newcommand{\unsur}[1]{\frac{1}{#1}}
\newcommand{\el}{\mathcal{L}}
\newcommand{\qq}{\mathcal{Q}}
\newcommand{\Per}{\mathrm{Per}}

\newcommand{\C}{\mathcal{C}}
\newcommand{\rond}{\hspace{-.1em}\circ\hspace{-.1em}}

\newcommand{\itm}{\item[-]}

\DeclareMathOperator{\supp}{Supp}

\DeclareMathOperator{\Int}{Int}

\newtheorem{prop} {Proposition} [section]
\newtheorem{thm}[prop] {Theorem} 
\newtheorem{defi}[prop] {Definition}
\newtheorem{lem}[prop] {Lemma}
\newtheorem{cor}[prop]{Corollary}
\newtheorem{theo}{Theorem} 
\newtheorem{coro}[theo]{Corollary}
\theoremstyle{remark}

\newtheorem{rmk}[prop]{Remark}

\begin{document}

\begin{abstract}
  We study the fine geometric structure of bifurcation currents in
  the parameter space of cubic polynomials viewed as dynamical
  systems. 
  In particular we prove that these currents
  have some laminar structure in a large region of parameter space,
  reflecting the possibility of quasiconformal deformations.  On the
  other hand, there is a natural bifurcation measure, supported on the
  closure of rigid parameters. We prove a strong non laminarity
  statement relative to this measure.
\end{abstract}

\maketitle

\section{Introduction}

The study of closed positive currents associated to bifurcations
 of  holomorphic dynamical systems is a topic under rapidly growing interest 
\cite{demarco, demarco2, bas-ber, pham, df}. The underlying thesis is
that bifurcation currents should contain a lot of information on the
geography of parameter spaces of rational maps of the Riemann sphere. 
 
The simplest  case  where the currents are not just measures
is when the parameter space under
consideration has complex dimension 2: a classical such example
is  the space $\mathrm{Poly}_3$ of cubic polynomials. The purpose of
the present article is to give a detailed account on the fine
structure of bifurcation currents and the bifurcation measure (to be defined 
shortly) in this setting. We chose to restrict
ourselves to cubics because, thanks to the existing literature, the
picture is much more complete in this case; nevertheless it is clear
that many results go through in a wider context.  

The work of Branner and Hubbard \cite{bh1, bh2} is a 
standard reference on the parameter space of cubic polynomials.  A crude way of
comparing our approach to theirs would be to say that in \cite{bh1,
  bh2} the focus was more on the topological properties of objects
while here we concentrate on their measurable and complex analytic
properties --this is exactly what positive closed currents do.

\medskip

Let us describe the setting more precisely (see Section \ref{sec:prel}
for details). Consider the cubic polynomial $f_{c,v}(z)=
z^3-3c^2z+2c^3+v$, with critical points at $\pm c$.  Since the
critical points play symmetric roles, we can focus  on $+c$. We
say that $c$ is {\em passive} near the parameter $(c_0, v_0)$ if the
family of holomorphic functions $(c,v)\mapsto f^n_{c,v}(c)$ is locally
normal. 
Otherwise, $c$ is said to be {\em active}. It is
classical that bifurcations always occur at parameters with active
critical points.

Let $\C^+$ be the set of parameters for which $+ c$ has bounded orbit
(and similarly for $\C^-$). It is a closed, unbounded subset in
parameter space, and it can easily be seen that the locus where $c$ is
active is $\fr\C^+$. We also denote by $\C= \C^+\cap \C^-$ the
connectedness locus.

The Green function $z\mapsto G_{f_{c,v}}(z)$ depends
plurisubharmonically on $(c,v)$ so the formulas $(c,v)\mapsto
G_{f_{c,v}}(\pm c)$ define plurisubharmonic functions in the space
$\cd_{c,v}$ of cubic polynomials. We define $T^\pm=dd^cG_{f_{c,v}}(\pm
c)$ as the {\em bifurcation currents} respectively associated to $\pm
c$. Notice that $ T_{\rm bif}= T^+ +T^-$ is the standard bifurcation
current, associated to the variation of the Lyapounov exponent of the
maximal entropy measure, as considered in 
\cite{demarco, demarco2, bas-ber, pham}. It is an easy fact that $\supp(T^+)$ is the
activity locus of $c$, that is, $\supp(T^+) = \fr\C^+$.

\medskip

The geometric intuition underlying positive   closed currents is that
of analytic subvarieties. Nevertheless, the ``geometry" of general
positive closed currents can be quite poor. 
{\em Laminar currents}  form a class of currents with rich geometric
structure, well suited for applications in complex dynamics. We say
that a positive closed current is {\em laminar} if it can be written
as an integral of compatible disks 
$T= \int [D_\alpha]d\nu(\alpha)$. Here compatible means that two disks
do not have isolated intersection points. The local geometry of a
laminar current can still be very complicated. 
On the other hand we say that $T$ is {\em
  uniformly laminar} if the disks $D_\alpha$ are organized as
laminations. The local geometry is tame in this case: a uniformly
laminar current is locally the ``product" of the leafwise integration
along the leaves by  a {\em transverse measure} on transversals. 

\medskip

The main point in our study is to look for some laminar structure for
the bifurcation currents in some regions of parameter space. To
understand why some laminarity should be expected, it is useful to
mention some basics on {\em deformations}. We say that two rational
maps are deformations of each other if there is a $J$-stable (in the sense of 
\cite{mss}) family
connecting them. There
is a natural stratification of the space of cubic polynomials,
according to the dimension of the space of deformations. The mappings
in the bifurcation locus correspond to dimensions 0 and 1. 
The dimension 1 case occurs 
for instance when one critical point is active and the other one is
attracted by an attracting cycle. The bifurcation locus contains a lot
of holomorphic disks near such a parameter, and it is not a surprise
that  the structure of the bifurcation current indeed reflects this fact. On the 
other hand we will see that the bifurcation measure $T^+\wedge T^-=T_{\rm bif}^2$
 concentrates on  the ``most bifurcating" 
part of the parameter space, that is, on the
 closure of the  set of parameters with zero dimensional deformation space. 

\medskip

Let us be more specific. Assume first that we are outside the
connectedness locus, say $-c$ escapes to infinity
(hence is passive). The structure of  $\C^+\setminus \C$ has been
described by Branner-Hubbard \cite{bh1, bh2}, and also Kiwi \cite{ki}
from a different point of view. There is a canonical
deformation of the maps outside $\C$, called the {\em wringing}
operation, and $\C^+$ locally looks like a holomorphic 
disk times a closed transversal set which is the union of countably many
copies of the Mandelbrot set and uncountably many points. 
Our first result is as follows (Proposition \ref{prop:lamin_escape}
and Theorem \ref{thm:dust}).

\begin{theo}\label{theo:outside}
The current $T^+$ is uniformly laminar outside the connectedness
locus, and the transverse measure gives full mass to the point
components. 
\end{theo}

The delicate part of the theorem is the statement about the
transverse measure. It follows from an argument of {\em similarity between
the dynamical and parameter spaces} at the measurable level, together
with a symbolic dynamics construction in the spirit of \cite{dds}. We also
give a new and natural  proof of the continuity of the  wringing operation, which is
a central result in \cite{bh1}. 

\medskip

Assume now that $-c$ is passive but does not escape to infinity (this
corresponds to parameters in $\mathrm{Int}(\C^-)$. This
happens when $-c$ is attracted towards an attracting cycle, and, if the
hyperbolicity conjecture holds, this is the only possible case.
We have the following result
(Theorem \ref{thm:dethelin} and Proposition \ref{prop:lamin_hyp}). 

\begin{theo}\label{theo:inside}
The current $T^+$ is  laminar in $\mathrm{Int}(\C^-)$, and uniformly
laminar in components where $-c$ is attracted to an attracting cycle.
\end{theo}

In particular laminarity holds regardless of the hyperbolic nature of
components. The proof is based on a general laminarity criterion 
due to De Th{\'e}lin \cite{dt}. The proof of uniform
laminarity in the hyperbolic case is more classical and 
follows from quasiconformal surgery. 

As a corollary of the two previous results, we thus get the following.

\begin{coro}
The current $T^+$ is laminar outside $\fr\C^+\cap \fr \C^-$.
\end{coro}

\medskip

In the last part of the paper, we concentrate on the remaining part of 
the parameter space. More precisely 
we study the structure of the {\em
bifurcation measure} $\mu_{\rm bif}=T^+\wedge T^-$, which is supported
in $\fr\C^+\cap \fr \C^-$. Understanding this
measure was already one of the main goals in \cite{df}, where it was
proved that $\supp(\mu_{\rm bif})$ is the closure of Misiurewicz
points (i.e. strictly critically preperiodic parameters),
and where the dynamical properties of $\mu_{\rm
 bif}$-a.e. parameters were investigated. It is an easy observation
 that $\mu_{\rm bif}$-a.e. polynomial has zero dimensional deformation space
(Proposition \ref{prop:rigidity}). 

It follows from \cite{df} that the currents $T^\pm$ are limits, in the
sense of currents, of the curves 
$$\Per^\pm(n,k)= \set{(c, v)\in \cd, \ f^n(\pm c)= f^k( \pm c)}.$$
 These curves intersect at Misiurewicz
points --at least the components where $\pm c$ are 
strictly preperiodic do. 
The laminarity of $T^\pm$ and the density of Misiurewicz
points could lead to the belief that the local structure of $\mu_{\rm
  bif}$ is that of ``geometric intersection" of the underlying
``laminations" of $T^+$ and $T^-$. The next theorem asserts that   
the situation is in fact more subtle (Theorem \ref{thm:isect}).  

\begin{theo}\label{theo:measure}
The measure $\mu_{\rm  bif}$ does not have local product structure on
any set of positive measure.
\end{theo}

Observe that this can also be interpreted as saying that generic
wringing curves do not admit continuations through $\fr \C$. 
As a perhaps surprising consequence of this result and previous work of
ours \cite{lamin, isect}, we obtain an asymptotic lower bound on the
(geometric) genus of the curves $\Per^\pm(n,k(n))$ as $n\cv \infty$, 
where $0\leq k(n)<n$ is any sequence. 
Recall that the geometric genus of a curve is the genus of its
desingularization; also, we define the genus of a reducible curve as
the sum of the genera of its components. Since these
curves are possibly very singular (for instance at infinity in the
projective plane), the genus cannot be directly read from the degree. 

\begin{theo} \label{theo:genus}
For any sequence $k(n)$ with $0\leq k(n)<n$, 
$$\unsur{3^n} \mathrm{genus} \left(\Per^\pm (n,k(n))\right)
\rightarrow \infty.$$
\end{theo}
  
\medskip

Finally, there is a striking analogy between the parameter space of
cubic polynomials with critical points marked
and dynamical spaces of polynomial automorphisms of
$\cd$. This has already be remarked  earlier:
 for instance, as far as the global topology 
of the space is concerned, the reader should
compare  the papers \cite{bh1} and \cite{ho}.  
We give in Table~\ref{table} a list of similarities and
dissimilarities between the two. Of course there is no rigorous link
between the two columns, and this table has to be understood more as a
guide for the intuition. For notation and basic concepts on polynomial automorphisms
of $\cd$, we refer to \cite{bs1, bls, fs, sib}.
On the other hand, Milnor's article \cite{mi} exhibits
similar structures between {\em parameter spaces} of cubic polynomials
and quadratic H{\'e}non maps.

\begin{table}\label{table}
\begin{footnotesize} 
\begin{center}
\begin{tabular}{|c|c|}
\hline
Parameter space of cubics & Polynomial automorphisms of $\cd$\\
 \hline\hline
$\C^\pm$, $\fr \C^\pm$, $\C$,  $\fr \C^+\cap \fr\C^-$ & $K^\pm$,
$J^\pm$, $K$, $J$\\
$T^+$, $T^-$, $\mu_{\rm bif}$ &    $T^+$, $T^-$, $\mu$\\
\hline
$\unsur{d^n}\left[\Per^\pm(n,k)\right]\cv T^\pm$ \cite{df}& 
$\displaystyle{\begin{cases}
\unsur{d^n}\left[(f^{\pm n})^*(L)\right]\cv T^\pm\\
 W^{s/u} \text{ are equidistributed}
\end{cases} }$ \cite{fs} \\
$\supp(\mu_{\rm bif})$ is the Shilov boundary of $\C$ \cite{df} & 
$\supp(\mu)$ is the Shilov boundary of $K$ \cite{bs1}\\
 $\supp(\mu_{\rm bif})= \overline{\set{\text{Misiurewicz pts}}}$ \cite{df}&   
$\supp(\mu)=\overline{\set{\text{saddle pts}}}$ \cite{bls}\\
$T^\pm$ are laminar outside  $\fr \C^+\cap \fr\C^-$ 
(Th.\ref{theo:inside} and \ref{theo:measure})& $T^\pm$ are
laminar \cite{bls}\\
Point components have full transverse measure (Th. \ref{theo:outside})
 &  The same when it makes sense \cite[Th. 7.1]{bs6}\\ 
 $T^+$ intersects all alg. subv. but the $\Per^+(n)$
(Prop. \ref{prop:alg_isect}) &
$T^+$ intersects all algebraic subvarieties \cite{bs2}\\
\hline
 $T^\pm$ are not extremal near infinity 
(Cor. \ref{cor:extremal}) & $T^\pm$ are extremal near infinity
\cite{fs, bs6}\\
$T^+\wedge T^-$ is not a geometric intersection
(Th.\ref{theo:measure}) &$T^+\wedge T^-$ is a 
 geometric intersection \cite{bls}\\
 \hline
 $\supp(\mu_{\rm bif})\subsetneq \fr \C^+\cap \fr\C^-$ & ??\\
Uniform laminarity outside $\C$ & ??\\
\hline
\end{tabular}
\end{center}
\end{footnotesize}
\caption{Correspondences between $\mathrm{Poly}_3$ and polynomial 
automorphisms of $\cd$.}
\end{table}
\medskip

We structure of the paper is as follows. In Section \ref{sec:prel} we
recall some basics on the dynamics on cubic polynomials and laminar currents. In
Section \ref{sec:escape}, we explain and reprove some results of
Branner-Hubbard and Kiwi on the escape locus, and prove Theorem
\ref{theo:outside}. Theorem \ref{theo:inside} is proved in Section
\ref{sec:lamin_bdry}. Lastly, in \ref{sec:nonlam} we discuss a notion
of {\em higher bifurcation} based on the dimension of the space of
deformations, and prove Theorems \ref{theo:measure} and \ref{theo:genus}.

\medskip

\noindent {\bf Acknowledgements.} This work is a sequel to \cite{df},
and would never have existed without the input of Charles Favre's
ideas. I also thank  Mattias Jonsson for pointing out to me  
the problem of understanding the bifurcation currents in the space of
cubic polynomials, and for several interesting discussions on the
topic. 


\section{Preliminaries}\label{sec:prel}

\subsection{Stability. Active and passive critical points}
The notion of stability considered in this article will be that of $J$-stability, 
in the sense of  \cite{mss}. 

\begin{defi}
Let $(f_\lambda)_{\lambda\in \Lambda}$ be a family of rational maps, parameterized by
 a complex manifold $\Lambda$. We say that $f_\lambda$ is $J$-stable (or simply
 stable) if there is a holomorphic motion of the Julia sets $J(f_\lambda)$, 
 compatible with the dynamics. 
 
Also, we say that $f_\lambda$ is a deformation of $f_{\lambda'}$ if there 
 exists a $J$-stable family connecting them.
\end{defi}

It is well known that if $f_\lambda$ is a deformation of $f_{\lambda'}$, then 
$f_\lambda$ and $f_{\lambda'}$ are quasiconformally conjugate in the neighborhoods
 of their respective Julia sets.
Notice that there is a stronger notion of stability, where the conjugacy is required on 
the whole Riemann sphere. 
This is the one considered for instance in \cite{mcms} for 
the definition of the Teichm{\"u}ller space of a rational map. The
stronger notion introduces some distinctions which are not relevant
from our point of view, like distinguishing the center from the other parameters in 
 a hyperbolic component.

\medskip

It is a central theme in our study to consider the bifurcations of critical points 
one at a time. This is formalized in the next classical definition. 

\begin{defi}\label{def:active}
Let $(f_\lambda, c(\lambda))_{\lambda\in \Lambda}$ 
be a holomorphic family of rational maps with a 
marked (i.e. holomorphically varying) critical point. 
 The marked critical point $c$ is \emph{passive} at $\lambda_0\in\Lambda$ if
  $\set{\lambda\mapsto f_\lambda^n c(\lambda)}_{n\in\nn}$ forms a
  normal family of holomorphic functions in the neighborhood of
  $\lambda_0$. Otherwise $c$ is said to be \emph{active} at
  $\lambda_0$. 
\end{defi} 

Notice that by definition the passivity locus is open while the activity locus is closed.
This notion is closely related to bifurcation theory 
 of rational maps, as the following classical 
proposition shows.

\begin{prop}[\cite{lyubactive, McM2}]\label{prop:mss}
Let $(f_\lambda)$ be a  family of rational maps with  all its critical points marked 
(which  is always the case, by possibly replacing 
$\Lambda$ with some branched cover). Then the family is $J$-stable 
iff all critical points are passive.
\end{prop}

\subsection{The space of cubic polynomials}\label{subs:param}

It is well known that the parameter space of cubic polynomials modulo affine
conjugacy has complex dimension 2. It has several possible presentations (see 
Milnor \cite{mi} for a more complete discussion). 

\medskip

The most commonly used parametrization is the following: for $(a,b)\in \cd$,
put \cite{bh1,mi}
$$f_{a,b}(z)= z^3-3a^2z+b.$$
The reason for the $a^2$ is that it allows the critical points ($+a$ and $-a$) 
to depend holomorphically on $f$.  

Two natural involutions in 
parameter space are of particular interest. 
First, the involution $(a,b)\mapsto (-a, b)$ preserves 
$f_{a,b}$ but exchanges the marking of critical points. In particular both critical points 
play the same role. 

The other interesting involution is  $(a,b)\mapsto (-a, -b)$. It sends $f_{a,b}$ to 
$f_{a,-b}$ which is conjugate to $f$ (via $-\mathrm{id}$). It can be shown that
the moduli space of cubic polynomials modulo affine conjugacy is 
$\cc^2_{a,b}/\left((a,b)\sim (-a, -b)\right)$, which is not a smooth manifold.

\medskip

In the present paper we will use the following parametrization (cf Kiwi \cite{ki}): 
$(c,v)\in \cd$, we put 
$$f_{c,v}(z)= z^3-3c^2z+2c^3+v = (z-c)^2(z+2c)+v.$$
The critical points here are $\pm c$, so they are both marked, 
and the critical values are 
$v$ and $v+4c^3$. The involution exchanging the marking of critical points 
then takes the form $(c,v)\mapsto (-c,  v+4c^3)$. 

Notice also that the two preimages of $v$ (resp. $v+4c^3$) are $c$ and
$-2c$ (resp. $-c$ and $2c$). The points 
$-2c$ and $2c$ are called {\it cocritical}.

The advantage of this presentation is that it is well behaved with respect to 
the compactification of  $\cd$ into the projective plane. 
More precisely, it separates the
sets $\C^+$ and $\C^-$ at infinity (see below Remark \ref{rmk:closure}).

\medskip

Notice finally that in \cite{df}, we used a still different parameterization, 
$$f_{c,\alpha}(z)= \unsur{3}z^3 - \unsur{2}cz^2+ \alpha^3,$$ where the
critical points are $0$ and $c$. In these  coordinates, the two
currents respectively associated to $0$ and $c$ have the same mass
(compare with Proposition \ref{prop:mass} below).

\subsection{Loci of interest in parameter space}\label{subs:loci}
We give a list of notation for the subsets in parameter space which will be of 
interest to us:
\begin{itemize}
\item[-] $\C^+$ the set of parameters for which $+c$ has bounded orbit. Similarly 
 $\C^-$ is associated to $-c$.
\item[-] $\Per^+(n)$ the set of parameters for which $+c$ has period $n$, and 
$\Per^+(n,k)$ the set of cubic polynomials
 $f$ for which $f^k(c)=f^n(c)$. Similarly for 
$\Per^-(n)$ and $\Per^-(n,k)$.
\item[-] $\C=\C^+\cap \C^-$ the connectedness locus.
\item[-] $\cd\setminus \C$ the escape locus.
\item[-] $ \cd\setminus(\C^+\cup\C^-)$ the shift locus.
\end{itemize}
 
The sets $\C^\pm$ and $\C$ are closed. Branner and Hubbard  proved in \cite{bh1}
 that the connectedness locus is compact and 
 connected\footnote{Our notation $\C^+$ corresponds to
  $\mathcal{B}^-$ in \cite{bh2}, and $\mathcal{E}^-$ in \cite{ki}.}.
 
 The $\Per^\pm(n,k)$ are  algebraic
curves. Since $\Per^+(n)\subset \C^+$, we see that $\C^+$ is
unbounded. 
It is easy to prove that the activity locus associated to $\pm c$ is
$\fr\C^\pm$. Consequently, the bifurcation locus is $\fr\C^+ \cup
\fr\C^-$.

\begin{rmk} Let $\om$ be a component of the passivity locus
associated to, say,  $+c$.  We say that $\om$ is hyperbolic if 
$c$ converges to an attracting cycle throughout $\om$. As an obvious
 consequence of the density of stability \cite{mss}, if the hyperbolicity conjecture holds, 
all passivity components are of this type. It would be interesting to describe their 
geometry. 
\end{rmk}


\subsection{Bifurcation currents}
\label{subs:currents}

 Here we apply the results of \cite{df} to define
 various plurisubharmonic functions and currents in parameter
 space, and list their first properties. We refer to Demailly's survey
 article \cite{de} for basics on positive closed currents. 
The  support of a current  or measure is denoted by $\supp(\cdot)$, and 
the trace measure by $\sigma_\cdot$.
Also $[V]$ denotes the integration current over
the subvariety $V$.

We identify the parameter
 $(c,v)\in \cd$ with the corresponding cubic polynomial $f$.
  For
 $(f,z)=(c,v,z)\in\cd\times\cc$, the 
 function $(f,z)\mapsto G_f(z)= G_{f_{(c,v)}}(z)$ is continuous and
 plurisubharmonic.  We define 
 \begin{equation} \label{eq:defg}
 G^+(c,v)=  G_{f}(c) \text{ and } G^-(c,v)=  G_{f}(-c). 
 \end{equation}
 The functions $G^+$ and $G^-$ are continuous, nonnegative, and plurisubharmonic.
 They have the additional property of being pluriharmonic when positive. 
 
 We may thus define (1,1) closed positive currents
 $T^\pm=dd^cG^\pm$, which will be referred to as the 
 {\em bifurcation currents associated to}
 $\pm c$. From \cite[\S 3]{df} we get that $\supp
 (T^\pm)$ is the activity locus associated to $\pm c$, that is,
  $\supp(T^\pm)=\fr\C^\pm$.
 
 It is also classical to  consider  the Lyapounov exponent
$ \chi= \log 3 + G^+ + G^- $
  of the equilibrium measure, and the
 corresponding current $T_{\rm bif}=dd^c\chi=T^++T^-$ known as the
 {\it bifurcation current} \cite{demarco, bas-ber}. Its support is the bifurcation locus.

\medskip

The following equidistribution theorem was proved in \cite{df}. 

\begin{thm}\label{thm:cv}
Let $(k(n))_{n\geq 0}$ be any sequence of integers such that  $0\leq k(n)<n$. Then 
$$\lim_{n\cv\infty} \unsur{3^n}\left[\Per^\pm(n,k(n))\right] = T^\pm.$$
\end{thm}

We can compute the mass of the currents $T^\pm$. The lack of symmetry
is not a surprise since  $c$ and $-c$ do not play the same role with respect to $v$. 

\begin{prop}\label{prop:mass}
The mass of $T^+$ with respect to the Fubini-Study metric on $\cd$ is
$1/3$. On the other hand the mass of $T^-$ is 1.
\end{prop}

\begin{proof} We use the equidistribution theorem. Indeed  a direct
  computation shows that the degree of $\Per^+(n)$ is $3^{n-1}$ while
  the degree of $\Per^-(n)$ is $3^{n}$. 
\end{proof}

We will mainly be interested in the fine geometric properties, in 
particular  {\em laminarity},  of the currents 
$T^\pm$. A basic motivation for this is the following observation.

\begin{prop}\label{prop:disk}
Let $\Delta$ be a holomorphic disk in $\cd$ where $G^+$ and $G^-$ are harmonic. 
Then the family  $\set{f\in \Delta}$ is $J$-stable.

This holds in particular when $\Delta\subset \C^+\setminus \fr\C^-$.
\end{prop}

\begin{proof} Recall that $G^+$ is nonnegative, and pluriharmonic where it is positive.
Hence if $G^+\rest{\Delta}$ is harmonic, either $G^+\equiv 0$  
(i.e. $c$ does not  escape) on $\Delta$, or $G^+> 0$  (i.e. $c$ always  escapes) on 
$\Delta$. In any case, $c$ is passive. 
Doing the same with $-c$ and applying Proposition \ref{prop:mss} then finishes the proof. 
\end{proof}

The laminarity results will  tell us that there are plenty of such disks, 
and  how these disks are organized in space:
 they will be organized as (pieces of) {\em laminations}. 

We can describe the relative positions
of $T^+$ and  algebraic curves, thus filling a line in Table \ref{table}.  Notice 
 that $G^+$ is continuous so $T^+\wedge [V]= dd^c(G^+\rest{V})$ is 
 always well defined. 
 Also, if $V\subset \supp (T^+)$, then $G^+\rest{V}=0$, so $T^+\wedge [V]= 0$. 
 
\begin{prop}\label{prop:alg_isect}
Let $V$ be an algebraic curve such that $T^+\wedge [V]=0$ then 
$V\subset \Per^+(n,k)$ for some $(n,k)$. If moreover 
$\supp(T^+)\cap V=\emptyset$, then $V\subset \Per^+(n)$ for some $n$.
\end{prop}

\begin{proof} If  $T^+\wedge [V]=0$, then $c$ is passive on $V$, which is an affine
algebraic curve,  
so by \cite[Theorem 2.5]{df}, $V\subset \Per(n,k)$ for some $(n,k)$. 
 
Assume now that  $\supp(T^+)\cap V=\emptyset$, and $V$ is a component of 
$\Per^+(n,k)$ where $c$ is strictly preperiodic at generic parameters (i.e. outside 
possibly finitely many exceptions). 
Then there exists  a persistent
cycle of length $n-k$ along $V$, on which $c$ falls after at most $k$ iterations. 
We prove that this leads to a contradiction.

The multiplier of this cycle defines a holomorphic function on $V$
which is thus constant.  If $f_0\in V$, $c$ is passive in the
neighborhood, so the multiplier cannot be greater than 1. On the other
hand if the multiplier is a constant of modulus $\leq 1$ along $V$, we
claim that the other critical point must have bounded orbit: indeed if
the cycle is attracting, there is a critical point in the immediate
basin of the cycle, which cannot be $c$ since $c$ is generically
strictly preperiodic. Now, if the multiplier is a root of unity or is
of Cremer type, then the cycle must be accumulated by an infinite
critical orbit, necessarily that of $-c$ (see e.g. \cite{mibook}), and
if there is a Siegel disk, its boundary must equally lie in the
closure of an infinite critical orbit. In any case, both critical
points have bounded orbits, so $K$ is connected. Since $V$ is
unbounded, this contradicts the compactness of the connectedness locus
\cite{bh1}.
\end{proof}


\subsection{Background on laminar currents}\label{subs:background}
Here we briefly introduce the notions of laminarity that will be
 considered  in the paper.
It is to be mentioned that the definitions of
flow boxes, laminations, laminar currents, etc. are tailored for the
specific needs of this  paper, hence not as general as they could be. 
 For more details on these notions, see  
\cite{bls, lamin, isect}; see also \cite{ghys} for general facts on
laminations and transverse measures. 

\medskip

We first recall the notion of direct integral of a family of positive closed currents.
Assume  that $(T_\alpha)_{\alpha\in \mathcal{A}}$ is a measurable family 
of positive closed currents in some open subset
$\om\subset\cd$, and $\nu$ is a
 positive measure 
on $\mathcal{A}$ such that (reducing $\om$ is necessary) $\alpha\mapsto \mathrm{Mass}(T_\alpha)$
is $\nu$-integrable. Then  we can define a positive closed current $T=\int T_\alpha 
d\nu(\alpha)$ by the obvious pairing with test forms
$$\langle T, \varphi\rangle = \int \langle T_\alpha, \varphi\rangle  d\nu(\alpha).$$

Let now $T$ and $S$ are positive closed currents in a ball $\om\subset\cd$. We say 
that the wedge product $T\wedge S$ is {\em admissible} if for some (hence for every) 
potential $u_T$ of $T$, $u_T$ is locally integrable with respect to the trace measure 
$\sigma_S$ of
$S$. In this case we may classically define $T\wedge S=dd^c(u_T S)$ which is a positive 
measure. The next lemma shows that integrating families of
positive closed currents  is well behaved with respect to taking wedge products.

\begin{lem}\label{lem:exchange}
Let  $T=\int T_\alpha d\nu(\alpha)$ as above, 
and assume that $S$ is a positive closed current such 
that the wedge product $T\wedge S$ is admissible.
Then for $\nu$-a.e. $\alpha$, $T_\alpha\wedge S$ is admissible
and 
$$ T\wedge S = \int (T_\alpha\wedge S) d\nu (\alpha).$$
\end{lem}

\begin{proof}
The result is local so we may assume $\om$ is the unit
  ball. Assume for the 
  moment  that there exists 
 a measurable family of nonpositive 
 psh functions $u_\alpha$, with $dd^cu_\alpha=T_\alpha$
 and $\norm{u_\alpha}_{L^1(\om')}\leq C \mathrm{Mass}(T_\alpha)$, for every $\om'
 \Subset \om$. 
In particular we get that $\alpha \mapsto \norm{u_\alpha}_{L^1(\om')}$ is 
$\nu$-integrable, and the formula $u_T= \int u_\alpha d\nu(\alpha)$ defines
a non positive potential of $T$. 
Then since $u_T\in L^1(\sigma_S)$ and all potentials are non positive, we
 get that for $\nu$-a.e. $\alpha$, 
$u_\alpha\in L^1(\sigma_S)$, and $u_TS = \int (u_\alpha S) d\nu(\alpha)$
by Fubini's Theorem.

\medskip

It remains to prove our claim. Observe first that if we are able to find a potential
$u_\alpha$ of $T_\alpha$, with 
 $\norm{u_\alpha}_{L^1(\om')}\leq C \mathrm{Mass}(T_\alpha)$, then, by 
 slightly reducing $\om'$, we can get a non positive potential with 
 controlled norm by substracting a constant, since by the submean inequality,
 $\sup_{\om''} u_\alpha \leq C(\om'') \norm{u_\alpha}_{L^1(\om')}$. Another observation
 is that it is enough to prove it when $T_\alpha$ is smooth and use regularization. 

 The classical way of finding a potential for a positive closed current of bidegree (1,1)
 in a ball is to first use the usual Poincar{\'e} lemma for $d$, and then solve a 
  $\overline{\fr}$  equation for a (0,1) form (see \cite{de}). The Poincar{\'e} lemma certainly
  preserves the $L^1$ norm (up to constants)
   because it boils down to integrating forms along 
  paths. Solving $\overline\fr$ with $L^1$ control in a ball is much more delicate 
  but still possible (see for instance \cite[p.300]{range}). This concludes the proof.
 \end{proof}

By {\em flow box}, we mean a closed family of disjoint holomorphic
graphs in the bidisk. In other words, it is the total space of a
holomorphic motion of a closed subset in the unit disk, parameterized
by the unit disk. 
It is an obvious consequence of Hurwitz' Theorem
that the closure of any family of disjoint graphs in $\dd^2$ is a  flow box. 
Moreover, the holonomy map along these family of graphs is
automatically continuous. This follows for instance from the
$\Lambda$-lemma of \cite{mss}. 

Let $\el$ be a flow box, written as the union of a family of disjoint
graphs as $\el= \bigcup_{\alpha\in \tau} \Gamma_\alpha$, where $\tau=
\el\cap (\set{0}\times\dd)$ is the central transversal. To every
positive measure $\nu$ on $\tau$, there corresponds a positive closed
current in $\dd^2$, defined by the formula 
\begin{equation}\label{eq:ul}
T= \int_\tau [\Gamma_\alpha]d\nu(\alpha). 
\end{equation}

A {\em lamination} in $\om\subset \cd$ is a closed subset of $\om$ which is locally
biholomorphic to an open subset of a flow box
$\el= \bigcup \Gamma_\alpha$. A positive closed 
current supported on a lamination is said to be {\em uniformly
  laminar} and {\em subordinate} to the lamination
if it is locally of the form \eqref{eq:ul}. 
Not every lamination carries a uniformly laminar current. This is the
case if and only if there exists an {\em invariant transverse
  measure}, that is, a family of positive measures on transversals,
invariant under holonomy. Conversely, a uniformly laminar current
induces a natural measure on every transversal to the lamination. 

\medskip

We say that two holomorphic disks $D$ and $D'$ are {\em compatible} if 
$D\cap D'$ is either empty or open in the disk topology. A positive
current $T$ in $\om\subset\cd$ is {\em laminar} if there exists a
measured family ($\mathcal{A},\nu$) 
of  holomorphic disks $D_\alpha\subset\om$, such that for every pair
$(\alpha,\beta)$, $D_\alpha$ and  $D_\beta$ are compatible  and 
\begin{equation}\label{eq_rep}
T=\int_\mathcal{A} [D_\alpha]d\nu(\alpha).
\end{equation}
The difference with uniform laminarity is that there is no (even
locally) uniform lower bound on the size of the disks in
\eqref{eq_rep}. In particular there is no associated
lamination, and the notion is strictly weaker. 
Notice also that the integral representation
\eqref{eq_rep} does not prevent $T$ from being closed, because of
boundary cancellation.

Equivalently, a current is laminar if it is the limit of an increasing
sequence of uniformly laminar currents. More precisely, 
$T$ is laminar in $\om$ if there exists a sequence of open subsets
$\om^i\subset\om$, such that   
$\norm{T}(\fr\om^i)=0$, together with an increasing sequence of
currents $(T^i)_{i\geq 0}$, $T^i$ uniformly laminar in $\om^i$,
converging to $T$. 

\medskip

In the course of section \ref{sec:nonlam}, we will be led to consider
{\em woven} currents. The corresponding definitions will be given at that time.


\section{Laminarity outside the connectedness locus}
\label{sec:escape}

In this section we give a precise description of $T^+$ outside the
connectedness locus. Subsections \ref{subs:lamin_wring} and \ref{subs:kiwi} are 
rather of expository nature.
We first recall the {\em wringing} construction
of Branner-Hubbard, and how it leads to uniform laminarity. 
Then, we explain some results of Kiwi \cite{ki} on the geometry of $\C^+$ at infinity.
In \S \ref{subs:dust},
based on an argument of similarity between the dynamical and parameter
spaces and a theorem of \cite{dds}, 
we prove that the transverse measure induced by $T^+$ gives
full mass to the point components. The presentation is as
self-contained as possible, only the results of \cite{bh2} depending on the combinatorics 
of tableaux are not reproved. 

\subsection{Wringing and uniform laminarity.}
\label{subs:lamin_wring}
We  start by defining some analytic functions in part of the parameter
space, by analogy with the definition of the functions $G^\pm$. Recall
that for a cubic polynomial $f$, the B{\"o}ttcher coordinate $\varphi_f$
is a holomorphic function defined in the open neighborhood of infinity
$$U_f=\set{z\in \cc, G_f(z)>\max(G^+, G^-)},$$ and semiconjugating $f$ to $z^3$
there, i.e. $\varphi_f\rond f= (\varphi_f)^3$. Also $\varphi_f=z+O(1)$
at infinity.

Assume $-c$ is the fastest escaping critical point,
i.e. $G_f(c)<G_f(-c)$. The corresponding critical value $v+4c^3$ has
two distinct preimages $-c$ and $2c$ (because $c\neq 0$), and it
turns out that 
$$\varphi_f(2c):= \lim_{U_f\ni z\cv 2c} \varphi_f(z) $$ is
always well defined (whereas there would be an ambiguity in defining  
$\varphi_f(-c)$). We put $\varphi^-(f)=\varphi_f(2c)$, so that
$\varphi^-$ is a holomorphic function in the open subset
$\set{G^+<G^-}$ satisfying
$G^-(f)= \log\abs{\varphi^-(f)}$. 

\begin{rmk} \label{rmk:continue} 
We notice for further use that  is possible to
continue $(\varphi^-)^3$ to the larger subset $\set{G^+<3 G^-}$ by
evaluating $\varphi_f$ at the critival value, and similarly for
$(\varphi^-)^9$, etc.  So locally, it is possible to define branches 
of $\varphi^-$ in larger subsets of parameter space by considering
inverse powers. 
\end{rmk}

Following Branner and Hubbard \cite{bh1} (we reproduce Branner's
exposition \cite{bra})
we now define the basic operation of {\em wringing the complex
structure}. This provides a holomorphic 1-parameter deformation of a
map in the escape locus, which has the advantage of being independent of
the number and relative position of escaping critical points. This is
in contrast with the full deformation theory 
\cite{mcms} which is sensitive to critical orbit relations, etc. \\

Let $\mathbb{H}$ be the right half plane, endowed with  the following 
group structure (with $1$ as unit element)
$$u_1*u_2=(s_1+it_1)*(s_2+it_2) =  (s_1+it_1)s_2+ it_2.$$ For $u\in
\mathbb{H}$ we define the diffeomorphism
 $g_u: \cc\setminus\dd \cv \cc\setminus\dd \cv$ by $g_u(z)= z\abs{z}^{u-1}$.    
Notice that $g_u$ commutes with $z^3$, and $g_{u_1*u_2} =
g_{u_1}\rond g_{u_2}$, i.e. $u\mapsto g_u$ is a left action on
$\cc\setminus \dd$. Since $g_u$ commutes with $z^3$ we can define a
new invariant almost complex structure on $U_f = \set{G_f(z)>\max(G^+,
  G^-)}$ by replacing it with
$g_u^*(\sigma_0)$ in the B{\"o}ttcher coordinate.  Here $\sigma_0$ is
the standard complex structure on $\cc$. More precisely define
$\sigma_u$ 
as the unique almost complex structure on $\cc$, invariant
by $f$, and  such that 
$$\sigma_u =
\begin{cases} 
\varphi_f^* g_u^*\sigma_0 &\text{ on } U_f\\
\sigma_0 &\text{ on } K_f
\end{cases}.$$
Of course, viewed as a Beltrami coefficient, 
$\sigma_u$ depends holomorphically on $u$, and it can be straightened
by a  quasiconformal map $\phi_u$ such that
\begin{itemize}
\item[-] $\phi_u\rond f\rond\phi_u^{-1}$ is a monic centered polynomial,
\item[-] $g_u\rond \varphi_f \rond \phi_u^{-1}$ is tangent to identity at
infinity. 
\end{itemize}
Under these assumptions $\phi_u$ is unique and depends holomorphically
on $u$ \cite[\S 6]{bh1}. \\

We put $w(f,u) = f_u= \phi_u\rond f\rond\phi_u^{-1}$. This is a
holomorphic disk trough $f=f_1$ in parameter space. Using the
terminology of \cite{dh}, $f_u$ is {\em hybrid equivalent} to $f$, so
by the uniqueness of the renormalization for maps with connected Julia
sets, if $f\in\C$, 
$f_u\equiv f$. On the other hand the construction
is non trivial in
the escape locus, as  the next proposition shows. 

\begin{prop} \label{prop:lamin_wring}
Wringing yields a lamination of $\cd\setminus\mathcal{E}$ by
Riemann surfaces, with leaves isomorphic to the disk or the punctured disk.
\end{prop}

This provides in particular
a new proof of the following core result of \cite{bh1}.

\begin{cor}\label{cor:lamin_wring}
The wring mapping $w:\cd\setminus\mathcal{E}
\rightarrow \cd\setminus\mathcal{E}$  is continuous.
\end{cor}

The proof of the proposition is based on the following simple but useful
lemma. 

\begin{lem}\label{lem:useful}
Let $\om\subset\cd$ be a bounded open set, and $\pi:\om\cv\cc$ be a 
holomorphic function. Assume that there exists a closed subset 
 $\Gamma$ and a positive constant $c$ so that 
\begin{itemize}
\item[-] For every $p\in \Gamma$, there exists a holomorphic map 
$\gamma_p:\dd\cv\Gamma$, with $\gamma_p(0)=p$ and $\abs{(\pi\rond
  \gamma_p)'(0)}\geq c$.
\item[-] the disks are compatible in the sense that $\gamma_p(\dd)\cap 
\gamma_q(\dd)$ is either empty or open in each of the two disks.
\end{itemize}
Then  
 $\overline{\bigcup_{p\in \Gamma}\gamma_p(\dd)}$ is a lamination  of
 $\Gamma$ by holomorphic curves. 
\end{lem} 

\begin{proof} We may assume $c=1$. First, since
 $(\pi\rond\gamma_p)'(0)$ never vanishes, $\pi$ is a
 submersion in the neighborhood of $\Gamma$. Hence  locally near $p\in
 \Gamma$, by choosing adapted coordinates $(z,w)$ we may suppose 
that $\pi$ is the first projection, i.e.  $\pi(z,w)=z$. 

The second claim is that $f$ is holomorphic and bounded by $R$ 
on $\dd$, with $f(0)=0$ and $\abs{f'(0)}\geq 1$, then there exists a
constant
 $\e(R)$ depending only on $R$ so that $f$ is injective on $D(0,
 \e(R))$.  Then by the Koebe Theorem, the image  univalently covers a disk of
radius $\e(R)/4$. 

Now observe that since $\om$ is bounded, if $\Gamma_1\subset\Gamma$ is relatively
compact, the family of mappings $\set{\gamma_p, \ p\in \Gamma_1}$ is
locally equicontinuous. Hence in the local coordinates $(z,w)$ defined
above, we can apply the previous observation to the family of
functions $\pi\rond \gamma_p$. The claim asserts that given any $z_0$,
and any point $p$ in the fiber $\Gamma\cap \pi^{-1}(z_0)$, the disk
$\gamma_p(\dd)$ contains a graph for the projection $\pi$
over the disk $D(z_0,\e/4)$, for some uniform $\e$. 

So all components
of   $\bigcup_{p\in \Gamma}\gamma_p(\dd)$ over the disk 
$D(z_0,\e/8)$ are graphs over that disk and they are disjoint. Moreover
they are uniformly bounded in the vertical direction by
equicontinuity. As observed in \S \ref{subs:background}, 
from Hurwitz'  Theorem and equicontinuity, we conclude that 
$\overline{\bigcup_{p\in \Gamma}\gamma_p(\dd)}$ is a lamination.   
\end{proof}

\begin{proof}[Proof of the proposition]
The idea is to use the previous lemma with $\varphi=\varphi^\pm$. We
need to understand how $\varphi^\pm$ evolve under wringing. We
  prove that the wringing disks form a lamination in
the open set $\set{G^+<G^-}$ where $\varphi^-$ is defined
--what we only need in the sequel is a neighborhood
of $\C^+$. The case of $\set{G^-<G^+}$ is of course
symmetric, and  remark \ref{rmk:continue} allows to extend the
reasoning to  a neighborhood of $\set{G^+=G^-}$. 

Recall that $w(f,u)=f_u=\phi_u\rond f\rond \phi_u^{-1}$, so if
$f=f_{(c,v)}$, the critical points of $f_u$ are $\phi_u(\pm
c)$. Also an explicit computation shows that 
$g_u\rond \varphi_f \rond \phi_u^{-1}$ semiconjugates $f_u$ and $z^3$ near
infinity, so by the normalization done, this is the B{\"o}ttcher
function $\varphi_{f_u}$. In particular, $\varphi_{f_u}(\phi_u(z))=
g_u(\varphi_f(z))$. 
 At the level of the Green function, we thus get
the identity
\begin{equation}\label{eq:green_wring}
G_{f_u}(\phi_u(z))= s G_f(z) \text{ where }u=s+it.
\end{equation}

Assume that $G^+<G^-$ so that the fastest escaping critical point is
$-c$. By (\ref{eq:green_wring}), the fastest  critical point of  $f_u$ is $\phi_u(-c)$ 
so 
\begin{equation}\label{eq:bott_wring}
\varphi^-(f_u)=g_u\rond \phi_f \rond \phi_u^{-1}(\phi_u(-c))= g_u(\varphi^-(f))=
\varphi^-(f)\abs{\varphi^-(f)}^{ u-1} 
\end{equation}
This equation tells us two things: first,
that $(f,u)\mapsto w(f,u)$ is locally uniformly bounded, and next, that 
 $f\mapsto \frac{\fr \varphi^-(w(f,u))}{\fr u}\rest{u=0}$ is locally uniformly bounded
from below in parameter space. 

Moreover the wringing disks are compatible because of the group action:
$w(w(f,v),u)=w(f,u*v)$. By using lemma \ref{lem:useful} we conclude 
that  they fit together in a lamination.\\

Now if $L$ is a leaf of the lamination, and $f\in L$ is any point,
$w(f,\cdot): \mathbb{H}\cv L$ is a universal cover.  Indeed it is
clearly onto because of the group action, and since both the leaf and
the map $w(f,\cdot)$ are uniformly transverse to the fibers of
$\varphi^-$, it is locally injective.
 
There are two possibilities:
either it is globally injective or not. In the first case $L\simeq
\dd$. In the second case, if $w(f,s_1+it_1)= w(f,s_2+it_2)$, then by 
(\ref{eq:green_wring}), $s_1=s_2$. Since near $w(f,s_1+it_1)$, the
leaf is a graph with uniform size for the projection $\varphi^-$, 
there exists a minimal $T>0$ such that $w(f,s_1+it_1+iT)= w(f,s_1+it_1)$. 
The group structure on $\mathbb{H}$ has  the property that $u*(v+iT)= 
(u*v) +iT$, from which we deduce that for every $u\in \mathbb{H}$ , $w(f, u+iT)
= w(f, u)$. In particular $L$ is conformally isomorphic to the
punctured disk. 
\end{proof}




\begin{proof}[Proof of Corollary \ref{cor:lamin_wring}] By definition,
the holonomy of a lamination is continuous.  We only have to worry about 
the compatibility between the natural parametrization 
of the leaves and the lamination
structure. Here, by (\ref{eq:bott_wring}),
the correspondence between the parameter $u$ on the 
leaves and the parameter induced by the transversal fibration 
$\varphi^{-}={\rm cst}$ is clearly uniformly bicontinuous --by uniform we mean here 
locally uniform with respect to the leaf. This proves the corollary.
\end{proof}

\begin{prop}\label{prop:lamin_escape}
The current $T^+$ is  uniformly laminar in the escape locus, and
subordinate to the  lamination by wringing curves.
\end{prop}

\begin{proof} The limit of a converging sequence of uniformly laminar
  currents, all subordinate to the same lamination is itself
 subordinate to this lamination: indeed locally in a flow box this is obvious. So by
 the convergence Theorem \ref{thm:cv}, to get the  uniform laminarity of $T^+$, 
it is enough to prove that the curves $[\Per^+(n)]$ are
 subordinate to the wringing lamination in
 $\cd\setminus\mathcal{E}$. But again this is obvious:  $w(f,u)$ is
 holomorphically conjugate to $f$ on $\Int(K_f)$, so  if $f$ has a
 superattracting cycle of some period, then so does $w(f,u)$. 
\end{proof}


\subsection{Geometry of $\C^+$ at infinity}\label{subs:kiwi}
In this paragraph we follow Kiwi \cite{ki}. For the sake of
 convenience, we include most proofs. 

The following lemma is \cite[Lemma 7.2]{ki}, whose proof relies on
easy explicit estimates. Besides notation, 
the only difference is the slightly more general hypothesis
$\abs{v}<3\abs{c}$.  We leave the reader check that this more general
assumption is enough to get the conclusions of the lemma.

\begin{lem} \label{lem:kiwi1}
If  $\abs{v}<3\abs{c}$ and $c$ is large enough, then the following hold:
\begin{itemize}
\itm $\abs{f^n(-c)}\geq \abs{c}^{3^n}(\sqrt{2})^{-3^{n-1}-1}$;
\itm $G^-(c,v)=\log\abs{c}+O(1)$;
\itm $G^+(c,v)\leq \unsur{3}\log\abs{c}+O(1)$.
\end{itemize}
\end{lem}

As a consequence we get an asymptotic expansion of 
$\varphi^-$ (compare \cite[Lemma 7.10]{ki}). Notice that under the
assumptions of the previous lemma, $G^+(c,v)<G^-(c,v)$ so $\varphi^-$
is well defined.

\begin{lem} \label{lem:kiwi2}
When $\abs{v}<3\abs{c}$ and $(c,v)\cv\infty$, then 
$$\varphi^-(c,v)= 2^{2/3}c + o(c).$$
\end{lem}

\begin{proof} Recall that 
$$\varphi^-(c,v)= \varphi_{f_{(c,v)}}(2c)= \lim_{n\cv\infty}
2c \left(\frac{f(2c)}{(2c)^3}\right)^\unsur{3}\cdots
\left(\frac{f^n(2c)}{(f^{n-1}(2c))^{3}}\right)^\unsur{3^n}.$$
For large $z$, $f(z)/z^3=1+O(1/z)$, hence since $f(2c)=f(-c)$, by
lemma \ref{lem:kiwi1} we get that for $n>1$, 
$$\frac{f^n(2c)}{(f^{n-1}(2c))^{3}}= 1+O\left(\frac{{2}^{3^n/2}}{\abs{c}^{3^n}}
\right),$$ and the product converges uniformly for large $c$. Moreover, 
${f(2c)}/{(2c)^3}$ converges to $1/2$ when 
$(c,v)\cv\infty$ and $\abs{v}\leq 3\abs{c}$, 
whereas, for $n\geq 2$,  ${f^n(2c)}/{(f^{n-1}(2c))^{3}}$ converges to 1. 
This gives the desired estimate.
\end{proof}

\begin{figure}
\begin{center}
 \psfrag{P1}{\small $[0:1:0]$}
 \psfrag{P2}{\small $[1:1:0]$}
 \psfrag{P3}{\small $[1:-2:0]$}
 \psfrag{L}{$L_\infty$}
 \psfrag{C}{$\mathcal{C}$}
 \psfrag{C+}{$\mathcal{C}^+$}
 \psfrag{C-}{$\mathcal{C}^-$}
 \psfrag{B}{${B}$}
 \psfrag{Phi}{\small $\varphi^-\!=k$}
\includegraphics[scale=0.7]{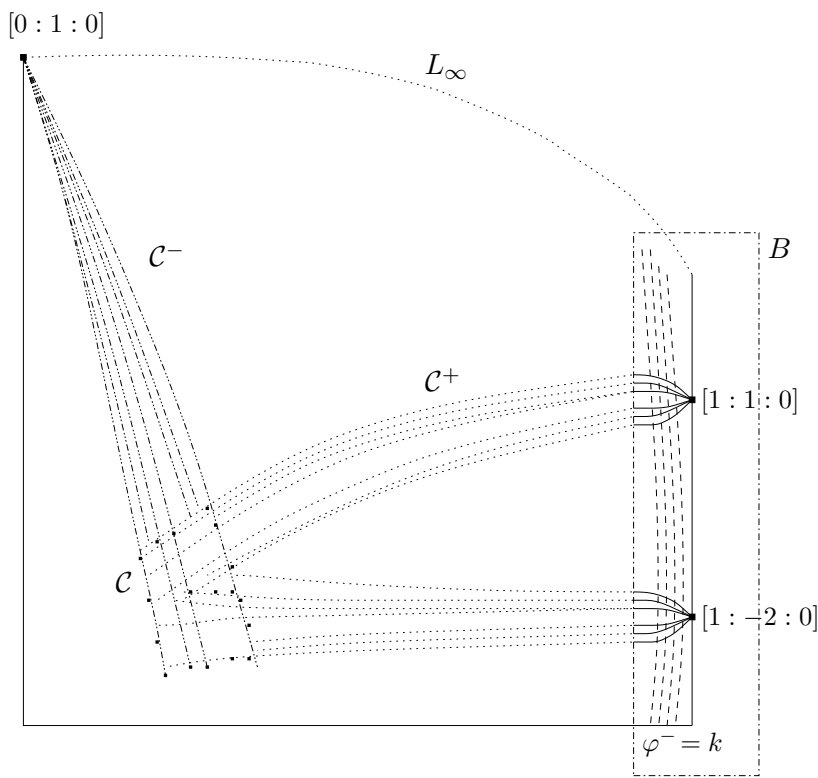}
\caption{Schematic picture of the parameter space.}\label{picture}
\end{center}
\end{figure}

We now translate these results into more geometric terms. See Figure
 \ref{picture} for a synthetic picture.
Compactify $\cd$ as
the projective plane $\pd$, and choose homogeneous coordinates 
$[C:V:T]$ such that $\set{[c:v:1], \ (c,v)\in\cd}$ is our 
parameter space. 
Consider the coordinates
$x=\frac{1}{c}=\frac{T}{C}$ and $y=\frac{v}{c}=\frac{V}{C}$ near
infinity, and  the bidisk 
$B=\set{\abs{x}\leq {x_0},\ \abs{y}\leq 1/3}$. 
The line at infinity becomes the vertical line $x=0$ in the new coordinates. 

\begin{prop}[Kiwi \cite{ki}]\label{prop:graphs}~
\begin{enumerate}[i.]
\item For ${k_0}$ large enough, the level sets $\set{\varphi^-=k, \ \abs{k}\geq 
k_0}$, are vertical holomorphic graphs in $B$. They fit 
together with the line at infinity $x=0$ as a holomorphic foliation. 
\item For small enough $x_0$, the leaves of the wringing lamination are graphs
over the first coordinate in  $B\setminus\set{x=0}$
\end{enumerate}
\end{prop}

\begin{proof} In the new coordinates,  
we have $$\varphi^-(x,y)= \frac{2^{2/3}}{x}(1+\delta(x,y)),$$ where $\delta(x,y)$ is a 
holomorphic function outside $\set{x=0}$. Since $\delta\cv 0$ as $x\cv 0$, $\delta$ 
extends holomorphically to $B$. Let $k'=1/k$; 
 the equation $\varphi^-=k$  rewrites as  
 $$x=2^{2/3}(1+\delta(x,y))k'.$$ For small enough $k'$ and fixed $y$, 
this equation  has exactly one solution in
 $x$, depending holomorphically on $y$, which means that 
 $\set{\varphi^-=k}$ is a vertical graph close to $x=0$. So clearly, 
$$\set{x=0}\cup \bigcup_{\abs{k}\geq k_0} \set{\varphi^-=k} $$ is a lamination near 
$\set{x=0}$. If we fix two small holomorphic transversals to $\set{x=0}$, the 
holonomy map is holomorphic outside the origin, so it extends holomorphically. 
We conclude that this lamination by vertical graphs
extends as  a holomorphic foliation across $\set{x=0}$.    \\
 
 On the other hand we know that the wringing curves are transverse to the 
 fibers of $\varphi^-$, i.e they are graphs for the projection 
$\varphi^-$. So by using the estimate of lemma  
\ref{lem:kiwi2} again and  Rouch{\'e}'s Theorem, we obtain that 
 the wringing leaves are graphs over the $c$ 
 coordinate for large $c$. Hence  the second part of 
the proposition. 
 \end{proof}

The next proposition asserts that the wringing curves contained in
$\C^+$ cluster only at two points at infinity. In particular, these
points are singular points for the wringing lamination. The results of
\cite[\S 7]{ki} may be thought as the construction of 
an abstract ``desingularization'' of this lamination.  

\begin{prop}[Kiwi \cite{ki}]\label{prop:closure}
The closure of $\C^+$ in $\pd$ is $\overline{\C^+} = \C^+\cup
\set{[1:1:0], [1:-2:0]}$. Likewise $\overline{\C^-}= \C^-\cup\set{
[0:1:0]}$. 
\end{prop}

\begin{rmk}\label{rmk:closure}
The choice of the $(c,v)$ parametrization is precisely due to this
proposition. By an easy computation, the reader may check
that in  $(a,b)$ coordinates of \S \ref{subs:param},  
$\C^+$ and $\C^-$ both cluster at the
same  point $[0:1:0]$ at infinity, an unwelcome feature. 
\end{rmk} 

\begin{proof} The starting point is the fact  that  $G^+ + G^-$ is
  proper in $\cd$: see \cite[\S 3]{bh1}.
 Consider now the following open neighborhood of $\C^+$:
$$\mathcal{V}=\set{ (c,v)\in\cd, \ G_f(v)=3G_f(c)<
  \unsur{3}G_f(-c)} = \set{3G^+<\unsur{3}G^-}.$$ We
will show that $\overline{\mathcal{V}}$ intersects the line at
infinity in $[1:1:0]$ and $ [1:-2:0]$, by proving that when
$G_f(-c)$ tends to infinity in $\mathcal{V}$, $v/c$ converges
to either $1$ or $-2$. The fact that both points are actually reached
is easy: this is already the case for the $\Per^+(2)$ curve.\\

Let us introduce some standard concepts. If $(c,v)\in \mathcal{V}$, 
in the dynamical plane, the level curves $\set{G_f(z)=r}$ are smooth
  Jordan curves for $r>G_f (-c)$, and $\set{G_f(z)=G_f (-c)}$ is a
  figure eight curve with self crossing at $-c$. By {\it disk at level} $n$,
  we mean a connected component of  
$\set{z,\ G_f(z)< 3^{-n+1}G_f (-c)}$. If $z$ is a point deeper than level
$n$, we define $D_f^n(z)$ as the  connected component of 
$\set{ G_f(z)< 3^{-n+1}G_f (-c)}$ containing $z$. Also 
$$A_f^0= \set{z,\ G_f (-c)< G_f(z)< 3G_f (-c)}$$ is an annulus of
modulus $$\unsur{2\pi} \log\frac{e^{3G_f (-c)}}{e^{G_f (-c)}}=
\unsur{\pi}G_f (-c).$$
If $(c,v)\in \mathcal{V}$, then $v$, $c$, and $-2c$ (the cocritical
point) are points of level $\geq 2$.
There are two disks of level 1 corresponding to the two inner components
of the figure eight: $D_f^1(c)$ and $D_f^1(-2c)$. 
In particular $D_f^1(v)\setminus\overline{D_f^2(v)}$ is an
annulus of modulus $\geq\unsur{2} {\rm modulus} (A_f^0)$ that
\begin{itemize}
\itm either separates $c$ and $v$ from $-2c$
\itm or separates $-2c$ and $v$ from $c$,
\end{itemize}
depending on which of the two disks contains $v$. 

When $(c,v)$ tends to infinity in $\mathcal{V}$, ${\rm modulus} (A_f^0)=G_f
(-c)\cv\infty $, and standard estimates in conformal geometry imply
that $\abs{v-c}=o(c)$ in the first case and $\abs{v+2c}=o(c)$ in the
second, that is,
$v/c$ respectively converges to $1$ or $-2$. To prove this, we  may
for instance use
the fact that an annulus with large modulus contains an essential
round annulus with almost the same modulus \cite{mcbook}.  \\

From this we easily get the corresponding statement for
$\mathcal{C}^-$, by noting that  the  involution
exchanging the markings
$(c,v)\mapsto (-c,v+4c^3)$ contracts the line at infinity with
$[1:0:0]$ deleted onto the point $[0:1:0]$. 
\end{proof}


\subsection{Transverse description of $T^+$}\label{subs:dust}

Branner and Hubbard gave in \cite{bh2} 
 a very detailed picture of 
$\C^+\setminus \C$ both from the point of view of the dynamics of an individual
mapping $f\in \C^+\setminus \C$ and the point of view of describing 
$\C^+\setminus \C$ as a subset in parameter space.  Roughly speaking, a map in 
$\C^+\setminus \C$ can be of two different types: quadratically renormalizable
or not. If $f\in \C^+\setminus \C$, $K_f$ is disconnected so it has infinitely many 
components. We denote by $C(+c)$ the connected component of the non escaping
critical point. Then
\begin{description}
\item[Renormalizable case] If  $C(+c)$ is periodic, then $f$ admits a
quadratic  renormalization, and 
$C(+c)$ is qc-homeomorphic to a quadratic Julia set. Moreover a component
of $K(f)$ is not a point if and only if it is a preimage of $C(+c)$.
\item[Non renormalizable case] If $C(+c)$ is not periodic, then $K_f$ is a Cantor set. 
\end{description}

\medskip

There is a very similar dichotomy in parameter space.  Notice that near infinity, $\C^+\subset B$, where $B$ is the
bidisk of proposition \ref{prop:graphs}. For large $k$, consider a disk
$\Delta$ of the form $\Delta=\set{\varphi^-=k} \cap B$. Then
$\mathcal{C}^+\cap\Delta$ is a disconnected compact subset of
$\Delta$. The connected components are of two different types
\begin{description}
\item[Point components]  corresponding to non quadratically renormalizable
maps.
\item[Copies of M] the quadratically renormalizable parameters are organized into
Mandel\-brot-like families \cite{dh}, giving rise to countably many
homeomorphic copies of the 
Mandelbrot set in $\Delta$.
\end{description}

\begin{rmk} \label{rmk:passive}
These results imply that  the (passivity) components
  of $\Int(\C^+)\setminus\C$ are exactly the open
  subsets obtained by moving holomorphically (under wringing) the
  components of the Mandelbrot copies . In particular, if the {\em  quadratic}
  hyperbolicity conjecture holds, all   these components are hyperbolic. 
\end{rmk}

 Since $T^+$ is uniformly laminar and
subordinate to the lamination by wringing curves, to have a good understanding 
of $T^+$, it is enough to
describe the transverse measure $T^+\wedge [\Delta]$. 

\begin{thm}\label{thm:dust}
The transverse measure induced by $T^+$ on a transversal 
gives full mass to the point components. 
\end{thm}

In particular, point components are dense in $\fr\C^+\cap \Delta$. 
Notice that the Mandelbrot copies are also dense, since they contain the points
$\Per^+(n)\cap\Delta$.

 \begin{proof}
 The proof  of the theorem will follow from similarity between the dynamical and 
 parameter spaces, and a symbolic dynamics argument in the style of \cite{dds}. 
 
In the domain of parameter space under consideration, $G^+<G^-$, so in
  the dynamical plane, the open set $\set{G_f(z)<G^-(f)}$ is bounded
  by a figure eight curve (see figure \ref{fig2}). 
Let $U_1$ and $U_2$ be its two connected
  components, and assume that $c\in U_2$, so that $f\rest{U_i}:U_i\cv
  \set{G_f<3G^-}$ is proper  with topological degree $i$. We denote by
$f_i$ the restriction $f\rest{U_i}$, $d_i=\mathrm{deg}(f_i)$ and by
  $U$ the topological disk $\set{G_f<3G^-}$.

\begin{figure}[h]
\begin{center}
  \psfrag{U1}{$U_1$}
  \psfrag{U2}{$U_2$}
  \psfrag{U}{$U$}
  \psfrag{+c}{\tiny$+c$}
  \psfrag{-c}{\tiny$-c$}
  \psfrag{fc}{\tiny$f(-c)$}
\includegraphics[scale=0.3]{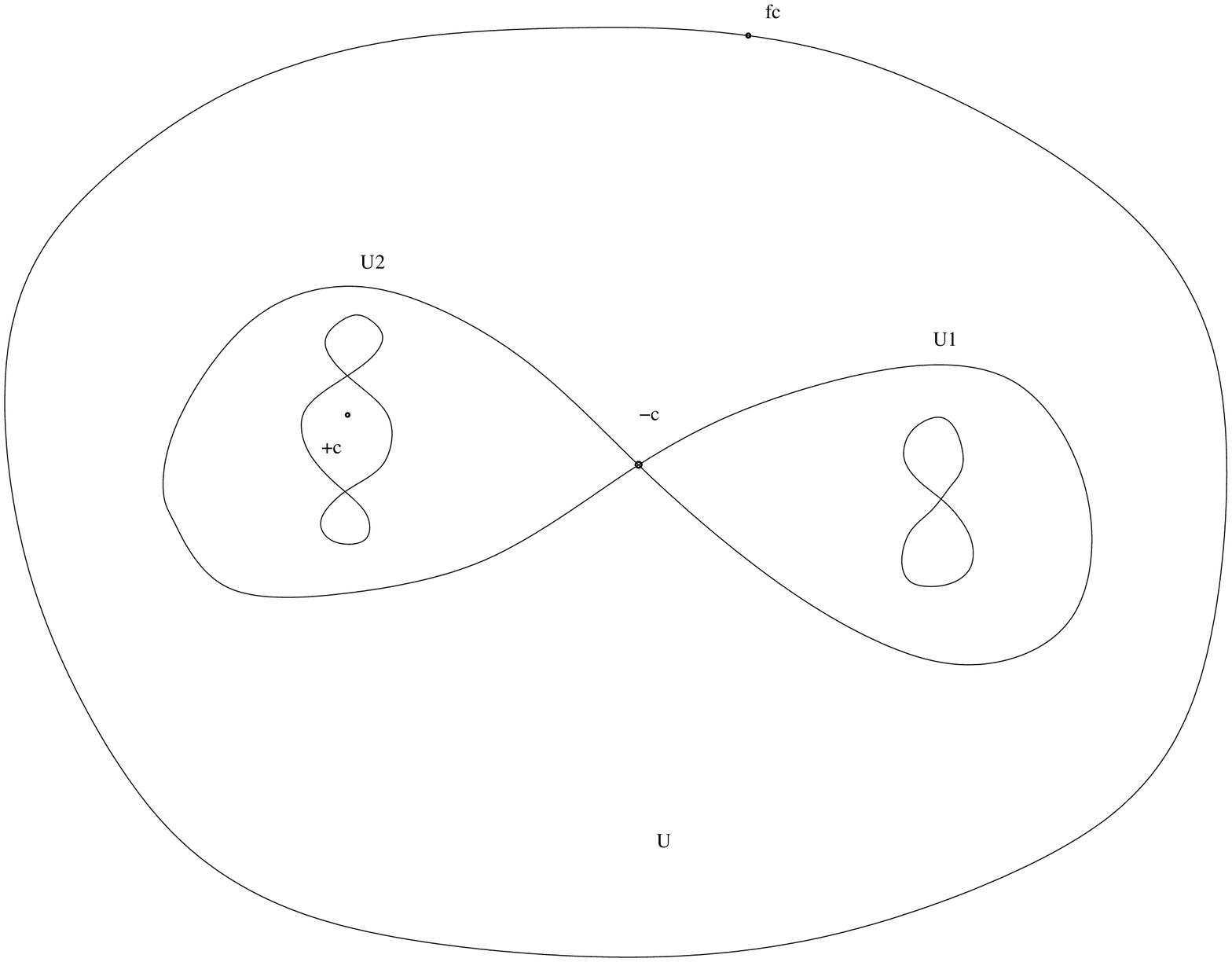}
\caption{}\label{fig2}
\end{center}
\end{figure}

Thus, classically, we get a decomposition of $K(f)$ in terms of
itineraries in the symbol space $\Sigma:=\set{1,2}^\mathbb{N}$. More
precisely,  
 for a sequence $\alpha\in\Sigma$, let 
$$K_\alpha= \set{z\in \cc,\ f^i(z)\in U_{\alpha(i)}}.$$ For every
$\alpha$, $K_\alpha$ is a nonempty closed subset, and the $K_\alpha$
form a partition of $K$. 

For example, for $\alpha=  \overline 2 = (2222\cdots)$, $K_{\overline
  2}$ is the filled Julia set of the quadratic-like map 
$f\rest{U_2}: U_2\cv U$. On the other hand $K_{\overline 1}$ is a
single repelling fixed point. Notice also that  $K_{\overline{12}}$
  is not a single point, even if the sequence $\overline{12}$ contains
  infinitely many 1's. 

\medskip

We now explain how the decomposition $K=\bigcup_{\alpha\in
  \Sigma}K_\alpha$ is reflected on the Brolin measure $\mu$. Let $p$
  be any point outside the filled Julia set and write 
\begin{equation}\label{eq:dds}
\unsur{3^n}(f^n)^*\delta_p= \unsur{3^n} 
\sum_{(\alpha_0, \ldots, \alpha_{n-1})\in
  \set{1,2}^n} \!\!\!  f_{\alpha_0}^*\cdots
  f_{\alpha_{n-1}}^* \delta_p.
\end{equation}
In \cite{dds}, we proved that for {\em any} sequence $\alpha\in \Sigma$,
and for any $p$ outside the Julia set,
the sequence of measures 
$$\unsur{d_{\alpha_0}\cdots d_{\alpha_{n-1}}} f_{\alpha_0}^*\cdots
  f_{\alpha_{n-1}}^* \delta_p $$
converges to a probability measure $\mu_\alpha$ independent of $p$,
  and supported on $\fr K_\alpha$. The measure $\mu_\alpha$ is the
  analogue of the Brolin measure for the sequence $(f_{\alpha_k})$.
This was stated  as Theorem
  4.1 and Corollary 4.6 in \cite{dds}
for sequences of horizontal-like in the unit bidisk but
  it is easy to translate it in the setting of polynomial-like maps in $\cc$.

 Let
  $\nu$ be the shift invariant measure on $\Sigma$, giving mass 
$(d_{\alpha_0}\cdots d_{\alpha_{n-1}})/3^n$ to the cylinder of
  sequences starting with $\alpha_0, \ldots, \alpha_{n-1}$ --the $(\unsur{3}, \frac23)$ 
  measure on $\Sigma$. Then
  \eqref{eq:dds} rewrites as 
$$\unsur{3^n}(f^n)^*\delta_p=\sum_{(\alpha_0, \ldots, \alpha_{n-1})\in
  \set{1,2}^n} \!\!\! \frac{d_{\alpha_0}\cdots
  d_{\alpha_{n-1}}}{3^n} \left[ 
\unsur{d_{\alpha_0}\cdots d_{\alpha_{n-1}}} f_{\alpha_0}^*\cdots
  f_{\alpha_{n-1}}^* \delta_p\right],$$ so at the limit we get the
  decomposition  $\mu = \int_\Sigma \mu_\alpha d\nu(\alpha)$. 

\medskip

 We fix a transversal disk $\Delta$ as in Proposition \ref{prop:graphs}.
  If $\Delta$ close enough to
  the line at infinity,  $G^+>0$ on $\fr\Delta$. Observe
   that the measure $T^+\wedge
\Delta$ induced by $T^+$ on $\Delta$ is the ``bifurcation current" of
the family of cubic polynomials parameterized by $\Delta$. We denote
it by $T^+\rest{\Delta}$.  

Following \cite[\S 3]{df} consider $\widehat\Delta=\Delta\times\cc$ and 
 the natural map $\hat f:
\widehat\Delta\cv\widehat\Delta $, 
defined by $\hat f(\lambda, z) = (\lambda, f_\lambda
(z))$. Consider any graph $\Gamma_\gamma=\set{ \ z=\gamma(\lambda)}$ 
over the first coordinate such that $\gamma(\lambda)$ escapes under
iteration by $f_\lambda$ for $\lambda\in \Delta$ (for instance
$\gamma(\lambda)= -c(\lambda)$). Then the sequence of currents
$\unsur{3^n} (\hat
f^*)^n [\Gamma]$ converges to a current $\widehat T$ such that 
$T^+ \rest{\Delta}= (\pi_1)_*(\widehat{T}\rest{\Gamma_{+c}})$.   Here $\pi_1:
 \widehat\Delta\cv \Delta$ is the natural projection. 

Now there is no choice involved in the labelling of $U_1$ and $U_2$ , so 
since  $G^+<G^-$ on $\Delta$, the decomposition  
$f_\lambda^{-1}U = U_1\cup U_2$ can be followed continuously
 throughout $\Delta$.
So with obvious notation, we get for $i=1,2$ a map 
$\hat{f}_i: \widehat{U}_i\cv \widehat U$.
This gives  rise to a coding of $\hat{f}$ orbits, in the same way as before. 
We thus get a decomposition of    
$\unsur{3^n} (\hat
f^*)^n [\Gamma]$ as 
$$\unsur{3^n} (\hat f^*)^n [\Gamma] = \unsur{3^n} 
\sum_{(\alpha_0, \ldots, \alpha_{n-1})\in
  \set{1,2}^n} \!\!\!  \hat f_{\alpha_0}^*\cdots
  \hat f_{\alpha_{n-1}}^* [\Gamma]. $$
The  convergence theorem of \cite{dds} implies that for any 
$\alpha\in \Sigma$, the sequence of currents
$$\widehat{T}_{\alpha, n}=
\unsur{d_{\alpha_0}\cdots d_{\alpha_{n-1}}}  \hat f_{\alpha_0}^*\cdots
  \hat f_{\alpha_{n-1}}^* [\Gamma]$$ converges to a current $\widehat
  T_\alpha$ in $\Delta\times \cc$. Indeed, 
 for every $\lambda\in \Delta$, 
the sequence of measures $
\unsur{d_{\alpha_0}\cdots d_{\alpha_{n-1}}} (f_\lambda)_{\alpha_0}^*\cdots
  (f_\lambda)_{\alpha_{n-1}}^* \delta_{\gamma(\lambda)}$ converges 
  by \cite{dds}. On the other hand the currents $\widehat{T}_{\alpha,
    n}$ are contained in a fixed vertically compact subset of
  $\Delta\times \cc$, and all cluster values of this sequence have the
  same slice measures on the vertical slices $\set{\lambda}\times
  \cc$. Currents with horizontal support 
being determined by their vertical slices,
  we conclude that the sequence converges (see \cite[\S 2]{dds} for basics
  on horizontal currents).  

\medskip

So again at the limit we get a decomposition of $\widehat{T}$ as an
integral of positive closed currents, 
$\widehat{T}= \int_\Sigma \widehat{T}_\alpha d\nu(\alpha)$, where
$\nu$ is the $(\unsur{3},\frac23)$ measure on $\Sigma$. 
Furthermore, for every $\alpha$, the wedge product 
$\widehat{T}_\alpha\wedge [\Gamma_c]$
 is well defined. This follows for instance from a classical transversality
 argument: recall that on $\fr \Delta$, $c$ escapes so $\supp
 \widehat{T}_\alpha \cap \Gamma_{c} \Subset \Gamma_{c}$. In
 particular any potential  $\widehat G_\alpha$  of
 $\widehat{T}_\alpha$ is locally integrable\footnote{
 The results in \cite{dds}  actually imply that $\widehat{T}_\alpha$ has
 continuous potential for $\nu$-a.e. $\alpha$} on $\Gamma_{c}$. 
Hence by taking the
wedge product with the graph $\Gamma_{c}$ and projecting down to
$\Delta$, we get a decomposition 
\begin{equation}\label{eq:decompos}
T^+\rest{\Delta}=  (\pi_1)_*(\widehat{T}\wedge [\Gamma_{c}])
= \int_\Sigma (\pi_1)_*(\widehat{T}_\alpha\wedge [\Gamma_{c}])
 d\nu(\alpha) = 
 \int_\Sigma T^+_\alpha\rest{\Delta}d\nu(\alpha).
 \end{equation} 
The second equality is justified by Lemma \ref{lem:exchange}.
 
 To fix the ideas, $T^+_\alpha\rest{\Delta}$
 should be understood as the contribution of the combinatorics
 $\alpha$ to the measure $T^+\rest{\Delta}$ in parameter space, and is supported on
 the set of parameters for which $c\in K_\alpha$. Notice also that since by definition 
 $c\in U_2$, for $T^+_\alpha\rest{\Delta}$ to be non zero, it is necessary that 
 $\alpha$ starts with the symbol 2.

\medskip
 
 Furthermore, since there are no choices involved in the labelling of
 the $U_i$, the decompositions $K=\bigcup K_\alpha$ and $\mu=\int
 \mu_\alpha d\nu(\alpha)$ are invariant under wringing. Hence the
 decomposition $T^+=\int T^+_\alpha d\nu(\alpha)$, which was defined
 locally, makes sense as a global decomposition of $T^+$ in
 $\cd\setminus \C$.
 
We summarize this discussion in the following
proposition.

\begin{prop}\label{prop:decompos}
There exists a decomposition of $T^+$ as an integral of positive
closed currents in $\cd\setminus \C$, 
$$T^+=\int_\Sigma T^+_\alpha d\nu(\alpha),$$ where $\Sigma =
\set{1,2}^\mathbb{N}$ and $\nu$ is the $(\frac13, \frac23)$-measure
on $\Sigma$. Moreover $T^+_\alpha$ is supported in the set of
parameters where $c$ has itinerary $\alpha$ with respect to the
natural decomposition  $\set{G_f(z)<G^-}=U_1\cup U_2$, with $c\in U_2$.
\end{prop} 

We now conclude the proof of Theorem \ref{thm:dust}. 
We combine the previous discussion with the
results of \cite{bh2}. In the dynamical space, if a component of the
Julia set is not a point, then it is preperiodic. Hence its associated
sequence $\alpha\in \Sigma$ is preperiodic. Since the measure $\nu$
does not charge points, the set of preperiodic symbolic sequences has
measure zero. Consequently, by the decomposition $\mu = \int
\mu_\alpha d\nu (\alpha)$, the set of preperiodic components has zero
$\mu$-measure. 

In parameter space, if the connected component of $\lambda$ in
$\C^+\cap\Delta$  is
not a point then  $C(+c)$ is periodic. Hence the itinerary sequence of
$c$ is periodic, so once again this corresponds to a set of sequences
of zero $\nu$-measure. In view of the decomposition 
$T^+\rest{\Delta}=\int T^+_\alpha\rest{\Delta} d\nu(\alpha)$, 
we conclude that the set of such
parameters has zero $T^+\rest{\Delta}$-measure.
\end{proof}

The current $T^+$ is not extremal in a small neighborhood of the line at infinity 
because of Proposition \ref{prop:closure}.
Notice that the two branches of $T^+$ at $[1:1:0]$ (resp. $[1:-2:0]$)
correspond to sequences starting with the symbol 21 (resp. 22).
More generally we have the following corollary, which fills a line in Table \ref{table}. 

\begin{cor}\label{cor:extremal}
The current $T^+$ is not extremal any open subset of the escape locus
$\cd\setminus\C$.
\end{cor}  

\begin{proof} Since the measure $\nu$ on $\Sigma$ does not 
charge points, and the currents $T^+_\alpha$ have disjoint supports
the decomposition of Proposition \ref{prop:decompos}
is non trivial in any open subset of $\cd\setminus\C$. 
\end{proof}


\section{Laminarity at the boundary of the connectedness locus}

\label{sec:lamin_bdry}

We start by a general laminarity statement in $\Int(\C^-)$
which, together with proposition \ref{prop:lamin_escape},
completes the proof of the laminarity of $T^+$ outside $\fr\C^-$.
Notice that since $\supp(T^+)=\fr\C^+$, this is a really a statement
about the boundary of the connectedness locus. 

\begin{thm}\label{thm:dethelin}
The current $T^+$ is laminar in $\Int(\C^-)$. 
\end{thm} 

\begin{proof} Consider a parameter $f_0\in \Int(\C^-)\cap \supp(T^+)=
  \Int(\C^-)\cap \fr\C^+\subset\fr\C$. The critical point
 $+c$ is active at $f_0$, whereas $-c$ is passive. Let
  $U\ni f_0$ be a small ball where $-c$ remains passive.
There are nearby parameters where $+c$ escapes, and we find that
 parameters in $U\setminus \C$ have a stable quadratic renormalization
  (see remark \ref{rmk:passive}). 

Because $T^+$ has continuous potential $G^+\geq 0$, we can write 
$$T^+=dd^cG^+=\lim_{\e\cv 0} dd^c\max(G^+, \e)=\lim_{\e\cv 0}T^+_\e.$$
Moreover $G^+$ is pluriharmonic where it is positive, so  
$dd^c\max(G^+, \e)$ has laminar structure. We will study the topology
of the leaves in $U$ and conclude by using a result by de Th{\'e}lin
\cite{dt}.\\

In the open set $U$, $-c$ never escapes, so the B{\"o}ttcher function
$\varphi^+(f)=\varphi_f(-2c)$ is well defined, and
$G^+=\log\abs{\varphi^+}$.  The real hypersurface $\set{G^+=\e}$ is foliated by
the holomorphic curves  $\set{\varphi^+=\exp({\e+i\theta})}$, and it is well
known (see e.g. \cite{de})
that $dd^c\max(G^+, \e)$ is uniformly laminar and subordinate to
this foliation --the fact that both $\set{G^+=\e}$ and the leaves are non
singular will be a consequence of the wringing argument below.  

The point is that in $U$ the leaves of this foliation are closed and
biholomorphic to the disk. Indeed, recall that the wring deformation 
$f\mapsto w(f,u)$ is a homeomorphism in $\cd\setminus \C$ for fixed
$u$: it can be inverted by using the group action. In particular
$w(\cdot,s)$ maps $U\cap \set{G^+=\e}$ homeomorphically 
into a neighborhood of infinity in $\Int(\C^-)$, where the foliation 
$\varphi^+=cst$ is well understood. By using proposition
\ref{prop:graphs} {\em i}, with $+$ and $-$ swapped, we get that the leaves
are planar Riemann surfaces. We conclude that the leaves are disks. Indeed
if $L=\set{\varphi^+=\exp({\e+i\theta})}$ is such a leaf in $U$,
$w(\cdot,s)(L)$ is an open subset in a vertical disk $\Delta$ (see
proposition \ref{prop:graphs} {\em i}), hence 
 $w(\cdot,s)^{-1}(\Delta)$ is a disk traced on
 $\set{\varphi^+= \exp({\e+i\theta})}$, containing $L$. So
 its intersection with the ball $B$ is a disk by the maximum principle.   \\

De Th{\'e}lin's Theorem \cite{dt} asserts that if $\Delta_n$ is a sequence
of (possibly disconnected)
submanifolds of zero genus in the unit ball, then any cluster value
of the sequence of currents $[\Delta_n]/{\rm Area(\Delta_n)}$ is a laminar
current. Pick a sequence $\e_n\cv 0$. 
In our situation, the approximating currents $T_{\e_n}^+$ are {\em
  integrals} of families of holomorphic disks. 
So, for every $n$, 
we first approximate $T_{\e_n}^+$ by a finite combination  of
leaves and denote the resulting approximating sequence of currents 
 by $(S_{n,j})_{j\geq 0}$. For every $j$, $S_{n,j}$ is a normalized finite sum   
of currents of integration along disks. Hence by choosing an appropriate
subsequence we can ensure that $S_{n,j(n)}\cv T^+$, and  applying
de Th{\'e}lin's result finishes the proof.
\end{proof}

Using quasiconformal deformations, we obtain a much more precise result in
hyperbolic (hence conjecturally all) components. Notice that due to Remark 
\ref{rmk:passive}, the quadratic hyperbolicity conjecture is enough to ensure uniform
laminarity of $T^+$ outside $\fr\C^+\cap \fr\C^-$ .

\begin{prop}\label{prop:lamin_hyp}
The current $T^+$ is uniformly laminar in hyperbolic components of
$\Int(\C^-)$. 
\end{prop}

\begin{proof} The proof will use a very standard quasiconformal
  surgical argument (cf. e.g. \cite[VIII.2]{cg}). 
  The major step is to prove
  that the foliation  by disks of the form $\set{\varphi^+=\e}$
  considered in the previous proposition does not degenerate in the 
  neighborhood of $\fr\C^+$. 

Before going into the details, we sketch the argument. 
Assume $f_0\in \fr\C^+\cap \om$, where $\om$ is a 
hyperbolic component of $\Int(\C^-)$. The
multiplier of the attracting cycle defines a natural holomorphic 
function in $\om$. Using quasiconformal surgery, we construct a
 transverse section of this function through $f_0$, which is a limit
 of sections of the form $\set{\varphi^+=\e}$, uniform in the sense of Lemma \ref{lem:useful}. Then by Lemma  \ref{lem:useful}, we conclude that the 
 foliation $\set{\varphi^+=\e}$ extends as a lamination to  $\fr\C^+$.
  Then, approximating $T^+$ by $T_\e^+$ as before implies
  that $T^+$ is uniformly laminar. \\

So let $\om$ be a hyperbolic component of $\Int(\C^-)$ and
$f_0\in \om\cap \overline{\mathcal{E}}$. The critical point $-c$ is attracted
by an attracting periodic point of period $m$ that persists through
$\om$.
Let $\rho$ be the multiplier of the cycle; $\rho$ is a
holomorphic function on $\om$. Let $U_0, \ldots, U_m=U_0$ be the cycle
of Fatou components containing the attracting cycle.
 Since $f_0$ is in the escape locus or at its boundary, $+c$ is not
 attracted by the attracting cycle., so its orbit does not enter 
 $U_0, \ldots, U_m$.  

The map $f_0^m:U_0\cv U_0$ is conjugate to a Blaschke product of degree
2:  there exists a conformal map $\varphi:U_0\cv\dd$
such that
$$ \varphi\rond f_0^m\rond \varphi^{-1} = z
\frac{z+\lambda_0}{1+\overline{\lambda_0} z}=B_{\lambda_0}(z),$$ where
$\lambda_0=\rho(f_0)$ is the multiplier of the cycle. For simplicity
we assume that $m=1$, see \cite{cg} for the adaptation to the general
case and more details. 

Choose a small  $\e$ and then $r$ so 
that $\abs{\lambda_0}<1-\e<r<1$, and for every $\lambda\in D(0,1-\e)$ we have
$$\overline{D(0,r)} \Subset B_\lambda^{-1}(\overline{D(0,r)}).$$
For $\lambda\in D(0,1-\e)$, choose a smoothly varying family of diffeomorphisms 
$$\psi_\lambda: f^{-1}\varphi^{-1} (\overline{D(0,r)})\longrightarrow
B_{\lambda}^{-1} (\overline{D(0,r)})$$
with $\psi_{\lambda_0}=\varphi$
and such that 
\begin{itemize}
\item[-] $B_\lambda\circ \psi_\lambda = \varphi\circ f$ on 
$\fr f^{-1}\varphi^{-1}(\overline{D(0,r)})$, 
\item[-] $\psi_{\lambda}=\varphi$  in
 $\varphi^{-1}(\overline{D(0,r)})$.
\end{itemize}
We can thus define a map 
$g_\lambda:\cc\cv\cc$ by 
$$\begin{cases}
g_\lambda = \varphi^{-1}\circ B_\lambda \circ \psi_\lambda \text{ in
}  f^{-1}\varphi^{-1} (\overline{D(0,r)})\\
g_\lambda= f  \text{ outside }f^{-1}\varphi^{-1} (\overline{D(0,r)}
\end{cases}$$
By definition of $\psi_\lambda$, the two definitions match on 
$\fr f^{-1}\varphi^{-1}(\overline{D(0,r)})$. 

Now define an almost complex structure $\sigma_\lambda$
by declaring that $\sigma_\lambda= \psi_\lambda^*\sigma_0$ on 
$f^{-1}\varphi^{-1}(\overline{D(0,r)})$, extending it by invariance under
$g_\lambda$, and 
 $\sigma_\lambda=\sigma_0$ outside
 the grand orbit (under $f_0$) 
of the attracting Fatou component. By construction, any orbit
under $g_\lambda$ hits the region where $g_\lambda$ is not holomorphic
at most once. So 
we can  straighten $\sigma_\lambda$, we obtain that  $g_\lambda$
is quasiconformally conjugate on $\cc$ to a cubic polynomial
$f_\lambda$, which  varies continuously with $\lambda\in D(0,1-\e)$. 
By construction, $f_{\lambda_0}=f_0$ and $\rho(f_\lambda)=\lambda$
(because the conjugacy is holomorphic in the neighborhood of the
attracting fixed point). \\

The family of maps $f_\lambda$ that we have constructed defines 
 a continuous section of
$\rho$ through $f_0$. Now, if $f_0 \notin \C$,  $\varphi^+(f_\lambda)$
is constant because $f_0$ and $f_\lambda$ are holomorphically
conjugate outside the Julia set. So $f_\lambda$ is constrained to stay
in the one dimensional local
leaf $L=\set{\varphi^+=\varphi^+(f_0)}$. Since $\rho$ is holomorphic on
  $L$ and $f_\lambda$ is a section of $\rho$, we obtain 
   that $f_\lambda$ depends holomorphically on $\lambda$, and that
   $\rho$ is a local biholomorphism on $L$. \\

 Now assume $f_0$ is at the boundary of the escape locus (that is at
 the boundary of $\C^+$). 
  Let  $(f_n)$ be a sequence of maps in $\mathcal{E}$ converging to
 $f_0$, say with $\rho(f_n)=\lambda_0$. Attached to each $f_n$ there is
 a holomorphic disk $\Delta_n=\bigcup_{\lambda\in D(0,1-\e)} f_n(\lambda)$ with
 $\rho(f_n(\lambda))=\lambda$.  Notice that any two such disks are equal or disjoint. 
 Let $\Delta_0$ be any normal limit of the sequence 
 of disks $(\Delta_n)$. This disk defines a section of $\rho$ through $f_0$ so 
 in particular we deduce that $\rho$ is submersive at $f_0$. 
 
  We claim that $\Delta_0$ does not depend on the sequence 
 $(\Delta_n)$. Indeed $\rho$ is a 
 submersion at $f_0$, hence  a fibration in the neighborhood of $f_0$, so the conclusion 
  follows from Hurwitz'  Theorem. 

By construction, this holomorphic disk is contained in $\fr \C^+$,
because along $\Delta_n$, $G^+$ is constant and tends to
zero as $n\cv\infty$. So both critical points are passive on
$\Delta_0$, and the dynamics along $\Delta_0$ is $J$-stable. 

\medskip

In conclusion we have constructed a family of holomorphic sections of
$\rho$ in $\om\cap (\mathcal{E}\cup\fr\C^+)$. In $\mathcal{E}$, the
sections are compatible in the sense of lemma \ref{lem:useful},
because they are 
subordinate to the foliation by $\set{\varphi^+=cst}$. As seen before, since
the sections on $\fr\C^+$ are obtained by taking limits, they are also
compatible by Hurwitz' Theorem. Hence lemma \ref{lem:useful} applies and tells us that these
disks laminate $\om\cap (\mathcal{E}\cup\fr\C^+)$. Hence we have
proved that  the natural foliation of $\om\cap\mathcal{E}$ admits a
continuation, as a lamination, to the boundary of $\mathcal{E}$.\\

It is now clear that $T^+$ is uniformly laminar in $\om$, because
$T^+$ can be written as $T^+=\lim T^+_\e=\lim dd^c \max(G^+, \e)$,
and, as we have seen in the proof of the previous theorem, 
$T^+_\e$ is subordinate to the natural foliation of $\om\cap
\mathcal{E}$. 
\end{proof}


\section{Rigidity at the boundary of the connectedness locus}
\label{sec:nonlam}

\subsection{Higher order bifurcations and the bifurcation measure}
So far we have seen that many parameters in the bifurcation locus 
admit  a one parameter family of deformations.  Due to laminarity, 
this holds for $T_{\rm bif}$-a.e. parameter outside $\fr\C^+\cap \fr\C^-$.
In this section we will concentrate on rigid parameters, that is, cubic 
polynomials with no deformations. 

Let $\mu_{\rm bif}$ be the positive measure defined by 
$\mu_{\rm bif}= T^+\wedge T^-$. An easy calculation shows that 
$T_{\rm bif}^2= (T^++T^-)^2= 2\mu_{\rm bif}$. 
One major topic in \cite{df} was to study the properties of $\mu_{\rm bif}$, which is 
in many respects the right analogue for cubics of the harmonic measure of the
Mandelbrot set. In particular it was proved that $\supp (\mu_{\rm bif})$ is the closure of
Misiurewicz points, and is a proper subset of $\fr\C^+\cap \fr \C^-$. 

\medskip

We define the {\em second bifurcation locus} $\mathrm{Bif}_2$
 as the closure of the set of rigid parameters, 
that is, parameters that do not admit deformations. 
By the density of Misiurewicz points, it is clear that 
$\supp(\mu_{\rm bif}) \subset \mathrm{Bif}_2$. 
In the next proposition we get a considerably stronger 
statement. Notice that the result is valid for polynomials of all degrees.

\begin{prop}\label{prop:rigidity}
The set of parameters $f$ for which 
there exists a  holomorphic disk with  
  $f\in \Delta\subset \C$ has zero $\mu_{\rm bif}$-measure. 
\end{prop}

In view of Proposition \ref{prop:disk} this can be understood as a 
generic rigidity result.

\begin{cor}
A $\mu_{\rm bif}$ generic polynomial admits no deformations.
\end{cor}

\begin{proof}[Proof of Proposition \ref{prop:rigidity}]  
Recall from \cite[\S 6]{df} that $\mu_{\rm
  bif} = (dd^c G)^2$, where $G=\max(G^+, G^-)$. Now if $\Delta$ is a
  holomorphic disk contained in $\C$, $G=0$ on $\Delta$. It is
  known that if $E$ is a subset of $\supp (dd^cG)^2$ such that
  through every point in $E$ there is a holomorphic disk on which $G$
  is harmonic, then $E$ has zero $(dd^cG)^2$ measure (see
  \cite[Corollary A.10.3]{sib}.

\medskip

Another argument goes as follows: by \cite[Theorem 10]{df},
 $\mu_{\rm bif}$-a.e.  polynomial $f$ satisfies the Topological Collet
 Eckmann (TCE) property. Also both
 critical points are on the Julia set so the only possible
 deformations come from invariant line fields \cite{mcms}.  On the other hand
 the Julia set of a TCE polynomial has
 Hausdorff dimension strictly  less than 2 \cite{pr}
so it has no invariant line fields. 
\end{proof}

The next result gives  some rough information on $\mathrm{Bif}_2$. We do not 
know whether the first inclusion is an equality or not. 

\begin{prop}\label{prop:bif2}
The second bifurcation locus satisfies
$$\supp(\mu_{\rm bif})\subset \mathrm{Bif}_2 \subsetneq 
(\fr\C^+\cap \fr\C^-) \cup \mathcal{E}\subsetneq \fr \C,$$
where the set $\mathcal{E}$ is empty if the hyperbolicity conjecture holds.
\end{prop}

\begin{proof} 
The relationship of $\mathrm{Bif_2}$ 
with active and passive critical points is as follows:
\begin{itemize}
\item[-] If one critical point is passive while the other is passive,
  it is expected that the passive critical point will give rise to a
  modulus of deformation. This is the case when it is attracted by a
  (super-)attracting cycle. Otherwise, we do not know how to prove it,
  nevertheless it is clear that there are nearby parameters with
  deformations, and moreover, these will be generic in the measure
  theoretic sense (see Theorem \ref{thm:dethelin}).  The
  existence of such parameters contradicts the hyperbolicity
  conjecture anyway.
\item[-] When the two critical points are active, there may exist deformations. 
The list of possibilities is as follows (see \cite{mcms}):
\begin{enumerate}
\item There is a parabolic cycle attracting both critical points, and their grand orbits 
differ.
\item There is a Siegel disk containing a postcritical point. 
\item There is an invariant line field --this is of course conjectured not to happen.
\end{enumerate}
\end{itemize}

In \cite{df}, we gave an example (due to Douady) 
of  a cubic polynomial $f_0 \in (\fr\C^+\cap \fr\C^-)\setminus  \supp(\mu_{\rm bif})$. 
This map has the property of having a parabolic fixed point attracting both critical points, 
as in case (1) above. 
Furthermore, 
every nearby  parameter either has a parabolic point, and is in fact conjugate to 
$f_0$, or  has an attracting point. In particular 
locally there is a disk of such parameters, conjugate to $f_0$. It is 
then clear that such a parameter is not approximated by rigid ones, so 
$f_0\notin \mathrm{Bif}_2$.
\end{proof}

 We have no precise information on the location of parameters in case (2) above. 
Here  is a specific
question: let $\Per_1(\theta)$ be the subvariety of parameters with a fixed point of 
multiplier $e^{2i\pi\theta}$. Assume that $\theta$ satisfies a Diophantine condition 
so that every $f\in\Per_1(\theta)$
 has a Siegel disk. For some parameters (see Zakeri \cite{zakeri}), it happens that 
 one of the critical point falls into 
 the Siegel disk after finitely many iterations 
 (the other one is necessarily in the Julia set). 
 In parameter space, this gives rise
 to disks contained in $\fr\C^+\cap \fr\C^-$. Are these disks contained in 
 $\supp(\mu_{\rm bif})$?  In $\mathrm{Bif}_2$?

\medskip



\subsection{Woven currents}
In the same way as foliations can be generalized to webs, there is a
class of woven currents, extending laminar currents. They were
introduced by T.C. Dinh \cite{dinh}. 

\begin{defi}\label{def:ug}
  Let $\om\subset \cd$ be an open subset. A web in $\om$ is any family
  of analytic subsets of $\om$, with volume bounded by some constant
  $c$. 

 A positive closed current $T$ in $\om$ is uniformly woven iff
  there exists a constant $c$, and a web $\mathcal{W}$ as above, with leaves of
  volume bounded by $c$, and endowed
  with a positive measure $\nu$, such that
$$T = \int_\mathcal{W} [V] d\nu(V).$$
\end{defi}

Notice that by Bishop's Theorem, 
any family of subvarieties with uniformly bounded volume is
relatively compact for the Hausdorff topology on compact subsets of $\om$. 
In analogy with the laminar case, we also define general woven currents.

\begin{defi}
A current $T$ in $\om$ is said to be woven if there exists a
sequence of open subsets $\om_i\subset \om$ such that $T$ gives zero
mass to $\fr\om_i$, and an increasing sequence of currents
$(T_i)$, converging to $T$, and such that $T_i$ is uniformly woven in
$\om_i$.
\end{defi}


\subsection{Non geometric structure of $\mu_{\rm bif}$}
In this paragraph, we investigate the structure of the bifurcation measure 
$\mu_{\rm bif}$, and prove that it is not the ``geometric intersection" of the laminar 
currents $T^\pm$. As a consequence we will derive an asymptotic
 estimate on the genera of the $\Per(n,k)$ curves. 

It is  natural to wonder whether  generic
 wringing curves in $\fr \C^\pm$ admit  continuations across $\fr \C$
(compare with the work of  Willumsen \cite{willumsen}). 
Many wringing curves
are subordinate to algebraic subsets  so the continuation
indeed exists: this is the case for the $\Per^\pm (n,k)$ curves as
well as the subsets  in parameter space defined by the condition that 
 a periodic cycle of given length has a given indifferent multiplier. 
Similarly, we may wonder whether the disks constituting
 the laminar structure of $T^\pm$ on $\fr \C\setminus (\fr\C^+\cap
 \fr\C^-)$ admit  continuations across $\fr\C^+\cap
 \fr\C^-$.

We will look at the continuation property from the
point of view of the structure of the bifurcation measure. If $E$ is a
closed subset in $\supp (\mu)$, we say that
$\mu_{\rm bif}$ {\em has local product structure on} $E$ if there exists
uniformly laminar currents $S^\pm\leq T^\pm$ defined in a neighborhood
of $E$ such that $\mu_{\rm bif}\rest{E}=S^+\wedge S^- $. Abusing
conventions, we will extend this definition to the case where
$S^+$ and $S^-$ are merely uniformly woven currents  and 
$0<S^+\wedge S^-\leq \mu_{\rm bif}\rest{E}$. 

This
terminology is justified by the fact that the wedge 
product of uniformly laminar
currents coincides with the natural
geometric intersection. More specifically,  if $S^+$ and $S^-$ are
uniformly laminar currents in $\om$, written as integral of disks of
the form $$S^\pm = \int [\Delta^\pm]d\mu^\pm,$$ then 
if the wedge product $S^+\wedge S^-$ is well defined, 
$$S^+\wedge S^- = \int [\Delta^+\cap \Delta^-]d\mu^+d\mu^-,$$
where $[\Delta^+\cap \Delta^-]$ is the sum of point masses at isolated
intersections of $\Delta^+$ and $\Delta^-$. 
We refer the reader to \cite{bls, isect} for
more details on these notions. Analogous 
results hold for uniformly woven currents.

\medskip

The next theorem shows that the structure of the bifurcation measure
is somewhat more complicated. Loosely speaking, this means
 that, in the measure theoretic sense, 
 neither wringing curves nor the disks of Section \ref{sec:lamin_bdry}
admit continuations across $\supp(\mu_{\rm bif})$.
 
\begin{thm}\label{thm:isect}
The measure $\mu_{\rm bif}$ does not have local product structure on
any set of positive measure.
\end{thm}

\begin{proof} Assume $\mu_{\rm bif}$ has local product structure on a
  set $E$ of positive measure. Let $S^\pm= \int
  [\Delta^\pm]d\mu^\pm$ be uniformly woven
  currents as in the previous discussion, with $0 <S^+\wedge S^-\leq
  \mu_{\rm bif}\rest{E}$. It is classical (see e.g. \cite[Lemma 8.2]{bls}
  that the currents $S^\pm$ have continuous potentials. 
Here the $\Delta^\pm$ are possibly singular
  curves of bounded area. By lemma \ref{lem:exchange}, 
$S^+\wedge S^-$ admits a decomposition as $\int [\Delta^+]\wedge S^-
  d\mu^+$. Since each $\Delta^+$ has finitely many
  singular points, the induced measure $[\Delta^+]\wedge S^-$  gives
  full mass to the regular part of $\Delta^+$, so we get that the
  variety  $\Delta^+$ is smooth around $S^+\wedge S^-$-a.e. point. 
  
  In particular , by Proposition \ref{prop:rigidity} above, discarding
  a set of varieties of zero $\mu^+$ measure if necessary, we may
  assume that for every $\Delta^+$, $\Delta^+\setminus \C \neq
  \emptyset$. Outside $\C$, $T^+$ is uniformly laminar and subordinate
  to the lamination by wringing curves, so $S^+$ itself is subordinate
  to this lamination and the transverse measure of $S^+$ is dominated
  by that of $T^+$ (see \cite[\S 6]{bls}, \cite{approx}).

From this discussion and Theorem \ref{thm:dust} we
 conclude that for $\mu^+$-almost every $\Delta^+$,
$\Delta^+\setminus \C$ is contained in at most countably many wringing
 curves,  and the
polynomials in $\Delta^+\setminus \C$ have Cantor Julia sets. Removing
 a set of zero $\mu^+$-measure once again, we may assume
 that this property holds for all $\Delta^+$. 

\medskip

By \cite[Corollary 11]{df}, there exists a set of parameters of full 
 $\mu_{\rm bif}$-measure $E'\subset E$ for which 
the orbits of both critical points are dense in the Julia set (this is
a consequence of the TCE property). By Fubini's Theorem, we may
 further assume that every $\Delta^+$ intersects $E'$.

Now since any variety $\Delta^+$ is
contained in $\C^+$, the point $+c$ is passive on $\Delta^+$. By
\cite[Theorem 4]{df}, if there exists a parameter $f\in  
\Delta^+$ for which $+c$ is preperiodic, then 
\begin{itemize}
\item[-] it is persistently preperiodic
\item[-] it is persistently attracted by an attracting periodic point
  or lies in a persistent Siegel disk.
\end{itemize} 
The first property only occurs on countably many varieties while the
second contradicts the genericity assumptions made on $\Delta^+$. So
we infer that $+c$ is never preperiodic on $\Delta^+$, in which case
its orbit can be followed by a holomorphic motion. But on $E'$, the
orbit of $+c$ is dense on the Julia set, which is connected, while on 
$\Delta\setminus \C$ it is contained in a Cantor set. We have reached a
contradiction.
\end{proof}

Not much seems to be known about the geometry of the $\Per^\pm(n)$ and
$\Per^\pm(n,k)$ curves. An unpublished 
result of Milnor's \cite{mismooth} asserts that the 
$\Per^\pm(n)$ are smooth in $\cd$. 
For every $0\leq k<n$, $\Per(n,k)$ is an algebraic
curve of degree $3^n$. The geometry of such curves 
is unknown. In particular these curves might be very
singular (this can also happen for the $\Per^\pm(n)$ at  infinity) 
and have many irreducible
components, so their genus cannot be directly 
estimated. 

As a consequence of Theorem \ref{thm:isect} and previous work of ours
 \cite{lamin, isect}, we get a rough asymptotic estimate on the genus of these
 curves. Recall that the geometric genus of a compact singular Riemann
 surface is the genus of its desingularization. In the next theorem,
 genus means geometric genus and, as usual, $\Per^+$ may be replaced
 by $\Per^-$. 

\begin{thm}
Let $k(n)$ be any sequence satisfying $0\leq k(n)< n$, and write the
decomposition of $\Per^+(k(n),n)$ into irreducible components
as $ \Per^+(k(n),n) =\sum m_{i,n} C_{i,n}$. 
Then 
\begin{equation}\label{eq:genus}
\unsur{3^n} \sum_i m_{i,n}\cdot \mathrm{genus} (C_{i,n}) \cv \infty.
\end{equation}
Consequently  $\max_i \mathrm{genus} (C_{i,n})\cv \infty$. 
\end{thm}

\begin{proof} In \cite{lamin}, we proved that if  a sequence 
of algebraic curves $C_n$ in $\pd$ satisfies certain geometric estimates,
  then the  cluster values of the sequence of currents
  $[C_n]/\mathrm{deg}(C_n)$ are laminar.  In \cite{isect} we proved that
  if $T_1$ and $T_2$ are two such laminar currents with
  continuous potentials, then the wedge product measure $T_1\wedge
  T_2$ has local product structure. By Theorem \ref{thm:isect}, we
  know that $\mu_{\rm bif}$ does not have product structure  on any
  set of positive measure. Therefore, the geometric estimates of
  \cite{lamin} are not satisfied. By inspecting the geometry of the
  $\Per^+(n,k(n))$ curves, we will see that this leads to
  \eqref{eq:genus}. A delicate issue is that \cite{lamin} requires an assumption 
  on the singularities which may not be satisfied when $k(n)\neq 0$, so 
we need to generalize it slightly. This is were woven currents are needed.  
  
 We start with the simpler case where $k(n)=0$. Then, by
  Milnor's result, the curves $\Per^+(n)$ are smooth in $\cd$, so in
  $\pd$ they can have singular points only at $[1:1:0]$ and
  $[1:-2:0]$. 
Let  $\Per^+(n) =\sum m_{i,n} C_{i,n}$ be the decomposition into
  irreducible components. We have to count the number of "bad
  components" when projecting $C_{i,n}$ from a generic point in
  $\pd$. Remark that the $C_{i,n}$ do not intersect in $\cd$, so no
  additional bad components arise from intersecting branches of
  different irreducible components of $\Per^+(n)$. 
Also,  the number of local irreducible components at
  every singular point of $C_{i,n}$ is bounded by
  $\mathrm{deg}(C_{i,n})$, since a nearby line cannot have more than 
 $\mathrm{deg}(C_{i,n})$ intersection points with $C_{i,n}$. 
By  \cite[Prop. 3.3]{lamin}, the number of bad components for each
  $C_{i,n}$ is thus bounded by $4\mathrm{genus}(C_{i,n}) + 6
  \mathrm{deg}(C_{i,n})$. By summing this estimate for all
  irreducible components, we get that the number (with multiplicity)
  of bad components is
  bounded by $4\sum_i m_{i,n}\cdot \mathrm{genus} (C_{i,n}) + O(3^n)$,
  and by \cite[Prop. 3.4]{lamin}, if \eqref{eq:genus} does not
  hold for a certain subsequence $n_i\cv\infty$, we get that $T^+$ is
  laminar everywhere, and strongly approximable in the sense of
  \cite[Def. 4.1]{isect}. The same holds for $T^-$, since $+c$ and $-c$ are
  exchanged by a birational involution. In particular, by
  \cite{isect}, $T^+\wedge T^-$ must have local product structure,
  which contradicts Theorem  \ref{thm:isect}.

\medskip
 
The proof in the general case involves a modification of
  \cite{lamin, isect}, which was already considered in \cite{dt2} (see
  also \cite{dinh}). The required modification is contained in the
  following proposition. We define the genus of  a reducible Riemann surface
as the sum of the genera of its components.

\begin{prop} \label{prop:woven}
Let $C_n$ be a sequence of (singular) algebraic curves in
  $\pd$, with $\deg(C_n)\cv\infty$. Assume that $\mathrm{genus}(C_n)=O(\deg(C_n))$.
Then the cluster values of the sequence $[C_n]/\deg(C_n)$ are woven
currents. 
\end{prop} 

\begin{proof} This is just a careful
  inspection of \cite{lamin}. We fix a generic point in $\pd$, and
  consider the central projection
  $\pi_p:\pd\setminus\set{p}\cv\pu$. For a subdivision $\qq$ of $\pu$ by
  squares, and $Q\in\qq$, 
 we look for {\it good components}, that is, 
{\it irreducible} components of $C_n\cap \pi_p^{-1}(Q)$
  which are graphs over $Q$. The difference with \cite{lamin} is that 
  two intersecting graphs are now  considered to yield two good components
 here. For each
  irreducible component $C_{i,n}$ of $C_n$, the number of bad
  components is controlled by using the Riemann-Hurwitz formula
 for the natural map  $\widehat{C}_{i,n} \cv \pu$, where
  $\widehat{C}_{i,n}$ is the normalization of $C_{i,n}$. A
  straightforward adaptation of \cite[Prop. 3.3]{lamin}  shows that
  the number of bad components is bounded by
  $O(\mathrm{genus}(C_{i,n}))+O(\deg(C_{i,n}))$. 
By summing for all $C_{i,n}$ and noting
  that any union of good components is good, we get that under the
  assumption of the lemma, the total
  number of bad components is  $O(\deg(C_n))$. 

We can now proceed with the proof of \cite{lamin}, by replacing
``laminar" by ``woven" everywhere, and get that the cluster values are
woven currents. 
\end{proof}

The woven currents  we obtain in this way satisfy explicit estimates
which allow to adapt \cite{isect} and 
study their intersection. This was already observed 
 in \cite[\S 3]{dt2}
--be careful that woven currents are called geometric
there. 

\begin{prop}\label{prop:isect_woven}
Let $T_1$, $T_2$, be woven currents obtained as cluster values of
curves satisfying the assumption of Proposition
\ref{prop:woven}. Assume further that $T_1$ and $T_2$ have continuous
potentials in some open set $\om$. Then the wedge product $T_1\wedge
T_2$ has local product structure in $\om$.
\end{prop}

\begin{proof} The sheme is as follows:
 consider two distinct linear projections
$\pi_1$ and $\pi_2$ in $\cd$, and subdivisions by squares of size $r$ of
the projection bases. This gives rise to a subdivision by cubes of
size $O(r)$ in $\cd$, which we denote by $\qq$.  If $T$ is a woven
current, as given by Proposition \ref{prop:woven}, then for a generic
such subdivision $\qq$ and every $Q\in \qq$, there exists a uniformly
woven current $T_Q$ in $Q$ such that the mass of $T-\sum_{Q\in \qq}
T_Q$ is $O(r^2)$. If now $T_1$ and $T_2$ are two such currents, with
continuous potentials, 
it is then easy to adapt \cite[Theorem 4.1]{isect}
and get that the wedge product $T_1\wedge T_2$ is approximated by 
$\sum_{Q\in \qq} T_{1,Q}\wedge T_{2,Q}$, which has geometric
interpretation by Lemma \ref{lem:exchange}. We conclude that the measure
$T_1\wedge T_2$ has local product structure. 
\end{proof}

The conclusion if the proof of the theorem is now clear. If for some
sequence $\mathrm{Per}^+(n,k(n))$, the estimate \eqref{eq:genus} is
violated, then by the previous proposition, the measure $\mu=T^+\wedge
T^-$ would have product structure, which is not the case due to
Theorem \ref{thm:isect}.
\end{proof}


\end{document}